\documentclass{amsart}
\usepackage{amssymb,amsmath}
\newcommand{\RR}{{\mathbb R}}

\newcommand{\e}{\varepsilon}

\newcommand{\Lap}{\Delta}
\newcommand{\de}{\delta}
\newcommand{\del}{\partial}
\newcommand{\om}{\omega}
\newcommand{\Om}{\Omega}

\newcommand{{\loc}}{{\ell\mathrm oc}}

\newcommand{\calH}{{\mathcal H}}

\def\meanint{{\diagup\hskip -.42cm\int}}

\newtheorem{theorem}{Theorem}
\newtheorem{lemma}{Lemma}
\newtheorem{proposition}{Proposition}
\newtheorem{corollary}{Corollary}

\newtheorem{remark}{Remark}

\begin{document}

\title[Differentiability for the Neumann Problem]
{Differentiability of Solutions to the Neumann Problem with Low-Regularity Data via Dynamical Systems}
\author{Vladimir Maz'ya}
\address{Link\"oping University}

\author{ Robert McOwen}
\address{Northeastern University}
\date{February 15, 2016}

\maketitle

\begin{abstract}
We obtain conditions for the differentiability of weak solutions for a second-order uniformly elliptic equation in divergence form with a homogeneous co-normal boundary condition. The modulus of continuity for the coefficients is assumed to satisfy the square-Dini condition and the boundary is assumed to be differentiable with derivatives also having this modulus of continuity.
Additional conditions for the solution to be Lipschitz continuous or differentiable at a point on the boundary depend upon the stability of a dynamical system that is derived from the coefficients of the elliptic equation.


\end{abstract}

\addtocounter{section}{-1}
\section{Introduction}
\smallskip\noindent

 For $n\geq 2$, let $U$ be a Lipschitz domain in $\RR^n$ with  exterior unit normal $\nu$ on $\partial U$.
Given a point $p\in\partial U$, let $B$ be an open ball centered at $p$. We want to 
consider solutions of the uniformly elliptic equation in divergence form in $U\cap B$ with homogeneous co-normal boundary condition on $\partial U\cap B$:
\begin{equation}\label{eq:NeumannProblem}
\begin{aligned}
\del_i(a_{ij}\,\del_j u)=0 \quad\hbox{in}\  U\cap B, \\
\nu_i\,a_{ij}\,\partial_j u=0 \quad\hbox{on}\ \partial U\cap B.
\end{aligned}
\end{equation}
(Here and throughout this paper we use the summation convention on repeated indices.)
Let $C_{comp}^1(\overline{U}\cap B)=\{u\in C^1(\overline{U}\cap B): \hbox{supp\,$u$ is compact in $\overline{U}\cap B$}\}$.
Recall that a {\it weak solution} of (\ref{eq:NeumannProblem}) is a function $u\in H^{1,2}(U\cap B)$, i.e.\ $\nabla u$ is  square-integrable on $U\cap B$, that satisfies
\begin{equation}\label{eq:weakform}
\int_U  a_{ij}\,\partial_ju\, \partial_i\eta \,dx=0
\quad\hbox{for all}\ \eta\in C_{comp}^1(\overline{U}\cap B).
\end{equation}
However, for irregular coefficients $a_{ij}$, a weak solution of (\ref{eq:NeumannProblem}) need not have a well-defined normal derivative along $\partial U$, so the boundary condition in \eqref{eq:NeumannProblem} is not meaningful, and we must only work with the variational formulation \eqref{eq:weakform}. 
When the coefficients $a_{ij}$ are  bounded and measurable, the classical results of Stampacchia \cite{S} show that a solution of \eqref{eq:weakform} is H\"older continuous on $\overline{U}\cap B$.
We want to consider mild regularity conditions on $a_{ij}$ and the boundary $\partial U$ under which a  solution of \eqref{eq:weakform} must be Lipschitz continuous, or even differentiable, at a given point of $\overline{U}\cap B$. 

We shall assume that the modulus of continuity $\om$ for the coefficients satisfies the {\it square-Dini condition}
\begin{equation}
\int_0^1 \om^2(r)\,\frac{dr }{r}< \infty.
\label{eq:Sq-Dini}
\end{equation}
Under this condition, the regularity of weak solutions at interior points of $U\cap B$ was investigated in \cite{MM3}, and found to also depend upon the stability of a first-order dynamical system derived from the coefficients. In this paper, we shall investigate the regularity of weak solutions at points on $\partial U\cap B$ and find somewhat analogous results.
Without loss of generality, we may assume that the boundary point is the origin $x=0$, and by a change of independent variables we may arrange $a_{ij}(0)=\delta_{ij}$. It turns out that the conditions for differentiability at a boundary point depend upon both the coefficients $a_{ij}$ and the shape of the boundary $\partial U$ in a rather complicated way, so for the purposes of describing our results in this introduction, let us consider two special cases: {\bf I.} When the boundary is flat near $0$. \ {\bf II.} When the operator is just the Laplacian near $0$. 

\medskip\noindent
{\bf I.}
Since our results are local in nature, we may assume the domain is the halfspace 
  \[
  \RR^n_+=\{(\widetilde x,x_n): x_n>0\}=\{(x_1,\dots,x_{n-1},x_n): x_n>0\}.
  \]
We assume $u\in H^{1,2}_{\loc}(\overline{\RR^n_+})$,
  i.e.\ first-order derivatives are integrable over compact subsets of $\overline{\RR^n_+}$. For $x\in \overline{\RR^n_+}$, let us write $x=r\,\theta$ where $r=|x|$ and $\theta\in S_+^{n-1}=\{x\in\RR^{n}_+: |x|=1\}$.
  We shall find that the relevant first-order dynamical system is
  \begin{equation}
 \frac{d\varphi}{dt} + R(e^{-t})\,\varphi=0 \quad\hbox{for}\ T<t<\infty,
 \label{eq:DynSystm}
 \end{equation}
where  $R(r)$ is the $(n-1)\times(n-1)$ matrix  given by
\begin{equation}\label{def:R-halfspace}
\left[R(r)\right]_{\ell k}:= \meanint_{S_+^{n-1}}\left(a_{\ell k}(r\theta)-n\,\sum_{j=1}^n a_{\ell j}(r\theta)\,\theta_j\,\theta_k\right)\,ds_\theta
\quad\hbox{for}\ \ell,k=1,\dots,n-1.
 \end{equation}
 Here and throughout the paper, the slashed integral denotes mean value.
Following \cite{C}, we say that (\ref{eq:DynSystm}) is {\it uniformly stable} as $t\to \infty$ if for every $\e>0$ there exists a $\de=\de(\e)>0$ such that any solution $\phi$ of  (\ref{eq:DynSystm}) satisfying $|\phi(t_1)|< \de$ for some $t_1>0$ satisfies $|\phi(t)|<\e$ for all $t\geq t_1$. (Since (\ref{eq:DynSystm}) is linear, an equivalent condition for uniform stability may be formulated in terms of the fundamental matrix; cf.\ Remark \ref{rk:2} in Appendix D.)
Moreover, a solution of  (\ref{eq:DynSystm}) is {\it asymptotically constant} as $t\to\infty$ if there is a constant vector $\phi_\infty$ such that $\phi(t)\to\phi_\infty$ as $t\to\infty$. As discussed in \cite{MM3}, if $R(r)\,r^{-1}\in L^1(0,\e)$, then (\ref{eq:DynSystm}) is
both uniformly stable and all solutions are asymptotically constant, but in general these conditions may be independent of each other.
As we shall see in Theorem \ref{th:1} in Section 2: {\it if \eqref{eq:DynSystm}, \eqref{def:R-halfspace} is uniformly stable as $t\to\infty$, then every solution $u\in H_{\loc}^{1,2}(\RR^n_+)$ of \eqref{eq:weakform} with $U=\RR^n_+$ is Lipschitz continuous at $x=0$; if, in addition, every solution of the dynamical system 
\eqref{eq:DynSystm}, \eqref{def:R-halfspace} is asymptotically constant as $t\to\infty$, then $u$ is differentiable at $x=0$.} 
Examples show (cf.\ Section 4) that solutions of \eqref{eq:weakform} need not be Lipschitz continuous at $x=0$
if the dynamical system \eqref{eq:DynSystm}, \eqref{def:R-halfspace} is not uniformly stable as $t\to\infty$.

\medskip\noindent
{\bf II.}
We assume $a_{ij}=\delta_{ij}$ and $U$ is a Lipschitz domain whose curved boundary $\partial U$ contains $x=0$, and let $B$ denote a ball centered at $0$. By a rotation of the independent coordinates, we may assume that $\partial U$ is given near $x=0$ as the graph of a Lipschitz function $h$, i.e.\ 
$x_n=h(\widetilde x)$ where $h(\widetilde 0)=0$. Since Lipschitz functions are differentiable almost everywhere, its gradient $\widetilde \nabla h$ is well-defined. We need $\widetilde\nabla h$ to satisfy the condition that
\begin{equation}\label{h'-sqDini}
\sup_{|x|=r}|\widetilde\nabla h(\widetilde x)|\leq \om(r) \ \hbox{as $r\to 0$},
\end{equation}
where $\om(r)$ satisfies the square Dini condition \eqref{eq:Sq-Dini}.
We again require stability properties of the dynamical system \eqref{eq:DynSystm}, but now the $(n-1)\times(n-1)$ matrix $R(r)$ is given by 
\begin{equation}\label{R-curved-Lap}
 \left[R(r)\right]_{\ell k} =n\,\meanint_{S_+^{n-1}}  \frac{\partial h(r\theta)}{\partial x_\ell}\theta_n\theta_k\,ds_\theta.
 \end{equation}
As a special case of Theorem 2 in Section 3, we have:  {\it if \eqref{eq:DynSystm}, \eqref{R-curved-Lap} is uniformly stable as $t\to\infty$, then every  solution $u\in H^{1,2}(U\cap B)$ of \eqref{eq:weakform}  is Lipschitz continuous at $x=0$; if, in addition, every solution of the dynamical system 
\eqref{eq:DynSystm}, \eqref{R-curved-Lap} is asymptotically constant as $t\to\infty$, then $u$ is differentiable at $x=0$.} 

\medskip
It is possible to obtain analytic conditions at $p$ that imply the desired stability of \eqref{eq:DynSystm}; we can even obtain conditions under which a solution of \eqref{eq:weakform} must have a critical point, i.e.\ $\nabla u(p)=0$. For $n=2$, of course, 
\eqref{eq:DynSystm} is a scalar equation, so  conditions for uniform stability and solutions being asymptotically constant are easily obtained; this is done in Section 4. For $n>2$,  conditions may be obtained in terms of the largest eigenvalue $\mu(r)$ of the symmetric matrix $S(r)=-\frac{1}{2}(R(r)+R^t(r))$, where $R^t$ denotes the transpose of $R$. Let us mention two conditions on $\mu(r)$:
\begin{equation}\label{mu-cond1}
\int_{r_1}^{r_2} \mu(\rho)\frac{d\rho}{\rho} <K \quad\hbox{for all $0<r_1<r_2<\e$}
\end{equation}
and
\begin{equation}\label{mu-cond2}
\int_{r}^{\e} \mu(\rho)\frac{d\rho}{\rho} \to -\infty \quad\hbox{ as $r\to 0$.}
\end{equation}
As an application to {\bf I}, if $R$ is defined by \eqref{def:R-halfspace}, then
we show in Section 2 that: \eqref{mu-cond1} 
{\it implies that  every solution $u\in H_{\loc}^{1,2}(\RR^n_+)$ of \eqref{eq:weakform} with $U=\RR^n_+$ is Lipschitz continuous at $x=0$ and \eqref{mu-cond2}  implies that $u$ is differentiable at $x=0$ with  $\partial_j u(0)=0$ for $j=1,\dots,n$.}
As an application to {\bf II}, if $a_{ij}=\delta_{ij}$ and $R$ is defined by \eqref{R-curved-Lap}, then the results in Section 3 show that: \eqref{mu-cond1} 
{\it implies that  every solution $u\in H_{\loc}^{1,2}(U\cap B)$ of \eqref{eq:weakform} is Lipschitz continuous at $x=0$ and \eqref{mu-cond2}  implies that $u$ is differentiable at $x=0$ with  $\partial_j u(0)=0$ for $j=1,\dots,n$.}
Additional analytic conditions on the matrix $R$ itself that imply the desired stability of \eqref{eq:DynSystm} may be found in \cite{MM3}, but we shall not discuss them further since they apply in general to the dynamical system \eqref{eq:DynSystm} and are not peculiar to the Neumann problem that we consider here.

Now let us say something about the methods used to prove these results.  
First we note that the modulus of continuity $\om(r)$ is a continuous, nondecreasing function of $r$ near $r=0$, and we need to assume that $\om$ does not vanish as fast as $r$ when $r\to 0$, i.e.\ for some $\kappa>0$
 \begin{equation}\label{om-vanishing}
\om(r)r^{-1+\kappa}\ \hbox{is nonincreasing for $r$ near $0$.}
\end{equation}
 Our analysis of regularity at $0\in\partial\RR^n_+$ is analogous to the analysis in \cite{MM3} for an interior point, and we shall adopt similar notation to make the parallels clear. In particular,  we use a decomposition 
 \begin{subequations}\label{spectral-decomp}
\begin{equation}\label{spectral-sum}
u(x)=u_0(r)+\widetilde v(r)\cdot\widetilde x+w(x),
\end{equation}
where the scalar function $u_0$ and $(n-1)$-vector function $\widetilde v=(v_1,\dots,v_{n-1})$ are given by
\begin{equation}\label{u0,vk}
u_0(r):=\meanint_{S_+^{n-1}} u(r\theta)\,ds_\theta,
\qquad
{v}_k(r):=\frac{n}{r}\,\meanint_{S_+^{n-1}} u(r\theta)\,\theta_k\,ds_\theta
\quad\hbox{for}\ k=1,\dots,n-1.
\end{equation}
Note that the scalar function $w$ has zero mean and first moments on the half sphere:
\begin{equation}
\meanint_{S_+^{n-1}}w(r\theta)\,ds_\theta=0=\meanint_{S_+^{n-1}}w(r\theta)\,\theta_k\,ds_\theta\quad\hbox{for}\ k=1,\dots,n-1.
\label{eq:meanint-w}
\end{equation}
\end{subequations}
As we shall see, the assumption that the dynamical system \eqref{eq:DynSystm}, \eqref{def:R-halfspace} is uniformly stable as $t\to\infty$ not only implies
that $\widetilde v$ and $\widetilde x\cdot\widetilde v'$ are bounded as $r\to 0$, but that $|u_0(r)-u_0(0)|$ and $|w(x)|$ are both bounded by $r\,\om(r)$ as $r\to 0$. Thus we have
\[
u(x)=u(0)+\widetilde v(r)\cdot\widetilde x + O(r\,\om(r)) \quad\hbox{as $r\to 0$}
\]
with $\widetilde v(r)$ bounded,
which shows that $u$ is Lipschitz at $x=0$. If we also know that all solutions of \eqref{eq:DynSystm}, \eqref{def:R-halfspace}
are asymptotically constant as $t\to\infty$, then we shall show $\widetilde v(r)=\widetilde v(0)+o(1)$ as $r\to 0$, which proves that 
$u$ is differentiable at $x=0$.

It could be of interest to compare our results on the differentiability of solutions to the Neumann problem with  asymptotic expansions that have been obtained for solutions of the Dirichlet problem (cf.\ \cite{KM} and \cite{KM2} which more generally consider elliptic operators of order $2m$). It is important to observe that the analysis at a boundary point for the Neumann problem is more complicated than it is for the Dirichlet problem. The reason for this can be clearly seen in the case of $\RR^n_+$ for $n\geq 3$:  the dynamical system that controls the behavior of $\widetilde v$ in the decomposition \eqref{spectral-sum} is  $(n-1)$-dimensional, while the corresponding decomposition for the Dirichlet problem involves only the coefficient of $x_n$, and so leads to a scalar ODE.

Let us mention that the square-Dini condition has been encountered in a variety of contexts: the differentiability of functions \cite{SZ}, Littlewood-Paley estimates for parabolic equations \cite{FSW}, and the absolute continuity of elliptic measure and $L^2$-boundary conditions for the Dirichlet problem \cite{CL}, \cite{FJK}, \cite{F}, \cite{KP}.  In addition, let us observe that the  projection methods used here were not only used in \cite{MM3} but also in \cite{MM1} and \cite{MM2}.

\section{A Model Problem for the Laplacian in a Half-space}

In this section, we consider \eqref{eq:weakform} when the operator is the Laplacian and
 $U=\RR^n_+$. However, in order for these results to be useful in our study of variable coefficients, we need to introduce some inhomogenoue terms to our variational problem. We assume that  $\vec f,f_0\in L_{\loc}^p(\overline{\RR^n_+})$ for some $p>n$, i.e.\ $f$ is $L^p$-integrable over any compact set $K\subset \overline{\RR^n_+}$. For $p'=p/(p-1)$ let
 \[
 H_{\rm comp}^{1,p'}(\overline{\RR^n_+}):=\{\eta\in H^{1,p'}(\RR^n_+): \eta(x)=0 \ \hbox{for all sufficiently large $|x|$}\},
 \]
and define
\begin{subequations}\label{Neumann-Laplace}
\begin{equation}\label{def:F}
F[\eta]=\int_{\RR^n_+} (f_0\eta-\vec f\cdot\nabla \eta)\,dx
\quad\hbox{for }\ \eta\in H_{\rm comp}^{1,p'}(\overline{\RR^n_+}).
\end{equation}
We now want to find a solution $u\in H^{1,p}_{\loc}(\overline{\RR^n_+})$ of the variational problem
\begin{equation}\label{eq:variational}
\int_{\RR^n_+} \nabla u\cdot\nabla \eta \,dx + F[\eta]   =0
\quad\hbox{for all}\ \eta\in H_{\rm comp}^{1,p'}(\overline{\RR^n_+}).
\end{equation}
\end{subequations}
We can obtain the solution using the Neumann function $N(x,y)$, which is a fundamental solution for $\Delta$  satisfying $\partial N/\partial x_n=0$ for $x\in\partial\RR^n_+$ and $y\in \RR^n_+$. Using the method of reflection, it can be written  as
\begin{equation}\label{def:N}
N(x,y)=\Gamma(x-y)+\Gamma(x-y^*),
\end{equation}
where $\Gamma(x)$ is the standard fundamental solution for the Laplacian $\Delta$ and $y^*=(\tilde y,-y_n)$
is the reflection in the boundary of $y=(\tilde y,y_n)\in \RR^n_+$. Using $N(x,y)$, we obtain the solution of 
\eqref{Neumann-Laplace} as follows. First, replace $\eta(y)$ in \eqref{eq:variational}  by $\chi_R(y)N(x,y)$ (for fixed $x$), where $\chi_R(y)=\chi(|y|/R)$ with a smooth cutoff function $\chi(t)$ satisfying $\chi(t)=1$ for $t<1$ and $\chi(t)=0$ for $t>2$. This can be done since $\nabla N(x,y)=O(|x-y|^{1-n})$ as $|x-y|\to 0$ implies (for fixed $x$) we have
$\chi_R(y)N(x,y) \in H_{comp}^{1,q}(\overline{\RR^n_+})$ for all $q<n/(n-1)$. Since $p>n$ is equivalent to $p'<n/(n-1)$, we have $\chi_R(y)N(x,y) \in H_{comp}^{1,p'}(\overline{\RR^n_+})$ and, provided the functions $\vec f$ and $f_0$ decay sufficiently as $|x|\to\infty$, we can let $R\to\infty$ to obtain the following solution formula for the problem \eqref{Neumann-Laplace}:
\begin{equation}\label{Neumann-Laplace-solution}
u(x)=\int_{\RR^n_+} (N(x,y)f_0(y)-\nabla_y N(x,y)\cdot\vec f(y))\,dy.
\end{equation}
For example, \eqref{Neumann-Laplace-solution} is the solution for \eqref{Neumann-Laplace} if $f_0,\vec f$ have compact support in 
$\overline{\RR^n_+}$.

In fact, we shall require a further refinement of \eqref{Neumann-Laplace}, but first we need to discuss projections. 
For $g\in L^1_{\loc}(\overline{\RR^n_+}\backslash\{0\})$ and $r>0$, let $Pg(r,\theta)$ denote  the projection of $g(r\theta)$ onto the functions on $S_+^{n-1}$ spanned by $1,\theta_1,\dots,\theta_{n-1}$:
\begin{subequations}
\begin{equation}\label{def:P}
Pg(r,\theta):=\meanint_{S_+^{n-1}} g(r\phi)\,ds_\phi+n\sum_{m=1}^{n-1}\theta_m\, \meanint_{S_+^{n-1}} \,\phi_m \, g(r\phi)\,ds_\phi,
\end{equation}
where we have used
\begin{equation}\label{def:c_n}
\meanint_{S^{n-1}_+} \theta_m^2\,ds_\theta=\frac{1}{n} \quad\hbox{for}\ m=1,\dots,n.\footnote{To verify \eqref{def:c_n}, note that $\theta_1^2+\cdots+\theta_n^2=1$ implies $\int_{S^{n-1}}\theta_i^2 ds=|S^{n-1}|/n$ for $i=1,\dots,n$. In particular,
$\int_{S^{n-1}_+}\theta_n^2 ds=|S^{n-1}|/2n=|S^{n-1}_+|/n$, and for $i=1,\dots,n-1$, $|S^{n-1}_+|=(n-1)\int_{S^{n-1}_+}\theta_i^2\,ds+|S^{n-1}_+|/n$, which yields \eqref{def:c_n}.}
\end{equation}
\end{subequations}
Note  that  $P\,1=1$ and $P\,\theta_m=\theta_m$ for $m=1,\dots,n-1$.
For $k\geq 1$, if
$g\in C^k(\overline{\RR_+^n}\backslash\{0\})$ then we can easily check that 
$Pg\in C^k((0,\infty)\times S^{n-1}_+)$; moreover, we have $\partial (Pg)/\partial x_n=0$ on $\RR^{n-1}\backslash\{0\}$ since there is no $\theta_n$-term in the definition of $Pg$.
Let us summarize this last remark in the following:
\begin{lemma} \label{le:d_nPu=0}   
For $k\geq 1$, if $u\in C^k(\overline{\RR_+^n}\backslash\{0\})$ then $Pu\in C^k((0,\infty)\times S^{n-1}_+)$
and $\partial Pu/\partial x_n=0$ on $\RR^{n-1}\backslash\{0\}$.  
\end{lemma}

If  $g\in L^1_{\loc}(\overline{\RR^n_+})$ and  $f$ is a bounded function with compact support in $\overline{\RR^n_+}$, then $Pf$ also has compact support and hence the product $g\,Pf$ is integrable on $\RR^n_+$. In fact, it is easy to see by Fubini's theorem that
\begin{equation}\label{eq:int-gPf=int-fPg}
\int_{\RR^n_+} g\,Pf\,dx=\int_{\RR^n_+} f\,Pg\,dx.
\end{equation}
Of course, \eqref{eq:int-gPf=int-fPg} also holds if $g\in L^{p}_{\loc}(\overline{\RR^n_+})$ and
$f\in L^{p'}_{comp}(\overline{\RR^n_+})$. In particular, if $g,g_j\in L^p_{\loc}(\overline{\RR^n_+})$ satisfy $\int fg_j\,dx\to\int fg\,dx$ for every $f\in L^{p'}_{comp}(\overline{\RR^n_+})$, then $\int fPg_j\,dx\to\int fPg\,dx$ for every $f\in L^{p'}_{comp}(\overline{\RR^n_+})$. In fact, we claim more:
\begin{lemma} \label{le:Puj->Pu}   
 If $p>n$ and $g,g_j\in H^{1,p}_{\loc}(\overline{\RR^n_+})$ satisfy $ \int_{\RR^n_+}\vec f\cdot\nabla g_j \,dx \to  \int_{\RR^n_+}\vec f\cdot\nabla g \,dx$ for every $\vec f\in L^{p'}_{comp}(\overline{\RR^n_+})$, then $ \int_{\RR^n_+}\vec f\cdot\nabla Pg_j \,dx \to  \int_{\RR^n_+}\vec f\cdot\nabla Pg \,dx$ for every $\vec f\in L^{p'}_{comp}(\overline{\RR^n_+})$.
\end{lemma}
\noindent{\bf Proof.} Since $p>n$ we have $H^{1,p}_{\loc}(\overline{\RR^n_+})\subset C(\overline{\RR^n_+})$,
and by density we may assume $\vec f\in C^1_{comp}(\overline{\RR^n_+})$. 
We can integrate by parts, and apply the above argument with $f=\hbox{div}\vec f$ in $\RR^n_+$ and $f_n$ in $\RR^{n-1}$:
\[
\begin{aligned}
 \int_{\RR^n_+}\vec f\cdot\nabla P(g_j) \,dx &=-  \int_{\RR^n_+}\hbox{div}\vec f\; P(g_j) \,dx
 +\int_{\RR^{n-1}} (f_n\,P(g_j))|_{x_n=0}\,d\tilde x \\
 &= -  \int_{\RR^n_+} P(\hbox{div}\vec f)\; g_j \,dx
 +\int_{\RR^{n-1}} (P(f_n)\,g_j)|_{x_n=0}\,d\tilde x \\
 &\to  -  \int_{\RR^n_+} P(\hbox{div}\vec f)\; g \,dx
 +\int_{\RR^{n-1}} (P(f_n)\,g)|_{x_n=0}\,d\tilde x \\
&= -  \int_{\RR^n_+}\hbox{div}\vec f\; Pg \,dx
 +\int_{\RR^{n-1}} (f_n\,Pg)|_{x_n=0}\,d\tilde x \\
& =  \int_{\RR^n_+}\vec f\cdot\nabla Pg \,dx.  \qquad \Box
 \end{aligned}
\]

Note that \eqref{eq:int-gPf=int-fPg} also enables us to define $P$ on distributions; for example, for $F$ as in \eqref{def:F} we have $PF[\eta]=F[P\eta]$ for any $\eta\in H_{\rm comp}^{1,p'}(\overline{\RR^n_+})$. Now, for a function or distribution $g$, let us define
\begin{equation}
g^\perp = (I-P)g.
\end{equation}
In particular, for $F$ as in \eqref{def:F}, we can define the functional  $F^\perp$:
\begin{subequations}\label{eq:Neumann-Laplace-perp}
\begin{equation}\label{def:F-perp}
F^\perp[\eta]:=\int_{\RR^n_+} (f_0\,\eta^\perp-\vec f\cdot\nabla (\eta^\perp))\,dx
\quad\hbox{for }\ \eta\in H_{\rm comp}^{1,p'}(\overline{\RR^n_+}).
\end{equation}
Now we can state the required refinement of \eqref{Neumann-Laplace}: to find $w\in H^{1,p}_{\loc}(\overline{\RR^n_+})$ satisfying $Pw=0$ and
\begin{equation}\label{eq:variational-perp}
\int_{\RR^n_+} \nabla w\cdot\nabla \eta \,dx + F^\perp[\eta]  =0
\quad\hbox{for all}\ \eta\in H_{\rm comp}^{1,p'}(\overline{\RR^n_+}).
\end{equation}
\end{subequations}

We will need the projection $P$ of $N(x,y)$ with respect to $y$ (i.e.\ for fixed $x$). To compute this, we
first expand $N(x,y)$ in spherical harmonics $\{\tilde\varphi_{k,m}: m=1,\dots,\tilde N(k)\ \hbox{and}\ k=0,\dots\}$ on $S^{n-1}$, where $\tilde N(k)$ is the dimension of the space of spherical harmonics that are even in $x_n$. In fact, as we show in Appendix A, assuming $n\geq 3$ this yields 
\begin{subequations}\label{eq:N(x,y)}
\begin{equation}\label{eq:N-x<y}
\begin{aligned}
N(x,y)&=\frac{a_0}{|y|^{n-2}}+\frac{a_0(n-2)}{c_n}\frac{|x|}{|y|^{n-1}}\sum_{m=1}^{n-1}\hat x_m\hat y_m\,+ \\
&\sum_{k=2}^\infty \,\frac{|x|^k}{|y|^{n-2+k}}
\sum_{m=1}^{\tilde N(k)}a_{k,m}\,\tilde\varphi_{k,m}\left(\hat{x}\right)\,\tilde\varphi_{k,m}\left(\hat{y}\right) 
\quad \hbox{for $|x|<|y|$},
\end{aligned}
\end{equation}
and
\begin{equation}\label{eq:N-y<x}
\begin{aligned}
N(x,y)&=\frac{a_0}{|x|^{n-2}}+\frac{a_0(n-2)}{c_n}\frac{|y|}{|x|^{n-1}}\sum_{m=1}^{n-1} \hat x_m\hat y_m+ \\
&\sum_{k=2}^\infty \,\frac{|y|^k}{|x|^{n-2+k}}
\sum_{m=1}^{\tilde N(k)}a_{k,m}\,\tilde\varphi_{k,m}\left(\hat{x}\right)\,\tilde\varphi_{k,m}\left(\hat{y}\right) 
\ \hbox{for $|y|<|x|$}.
\end{aligned}
\end{equation}
\end{subequations}
The coefficients $a_0,a_{k,m}$\ can be computed but their values are not important to us now;
and we have here used the notation $\hat x=x/|x|$ and $\hat y=y/|y|$ (although elsewhere we have used $\theta=x/|x|$).
If we denote the projection $P$ of $N(x,y)$ with respect to $y$ simply by $PN(x,y)$, then we have
\begin{equation}\label{def:PN}
PN(x,y)=
\begin{cases}
\frac{a_0}{|y|^{n-2}}+\frac{a_0(n-2)}{c_n}\frac{|x|}{|y|^{n-1}}\sum_{m=1}^{n-1}\hat x_m\hat y_m \quad \hbox{for $|x|<|y|$},\\
\frac{a_0}{|x|^{n-2}}+\frac{a_0(n-2)}{c_n}\frac{|y|}{|x|^{n-1}}\sum_{m=1}^{n-1} \hat x_m\hat y_m  \quad \hbox{for $|y|<|x|$}.
\end{cases}
\end{equation}
We can also define
\begin{equation}
\begin{aligned}
N^\perp(x,y)&=N(x,y)-PN(x,y)\\
&=\begin{cases}
\sum_{k=2}^\infty \,\frac{|x|^k}{|y|^{n-2+k}}
\sum_{m=1}^{\tilde N(k)}a_{k,m}\,\tilde\varphi_{k,m}\left(\hat{x}\right)\,\tilde\varphi_{k,m}\left(\hat{y}\right) \quad \hbox{for $|x|<|y|$},\\
\sum_{k=2}^\infty \,\frac{|y|^k}{|x|^{n-2+k}}
\sum_{m=1}^{\tilde N(k)}a_{k,m}\,\tilde\varphi_{k,m}\left(\hat{x}\right)\,\tilde\varphi_{k,m}\left(\hat{y}\right) 
\ \hbox{for $|y|<|x|$}.
\end{cases}
\end{aligned}
\end{equation}
Using the same argument as for \eqref{Neumann-Laplace-solution}, provided the functions $\vec f$ and $f_0$ decay sufficiently as $|x|\to 0$ and $|x|\to\infty$, we have the following solution formula for the problem \eqref{eq:Neumann-Laplace-perp}:
\begin{equation}\label{eq:WeakSoln-Laplace-perp}
w(x)=\int_{\RR^n_+}  \left(N^\perp(x,y) f_0(y) - \nabla_y N^\perp(x,y) \cdot  \vec f(y) \right) \, dy.
\end{equation}
For example, \eqref{eq:WeakSoln-Laplace-perp}
holds if $f_0,\vec f$ have compact support in 
$\overline{\RR^n_+}$.

Let us now obtain estimates on the solution of \eqref{eq:Neumann-Laplace-perp} given by \eqref{eq:WeakSoln-Laplace-perp} when we make certain assumptions about the decay of $f_0$ and $\vec f$ as $|x|\to 0$ and $|x|\to \infty$. We do so using the $L^p$-mean on annuli: for  $r>0$ define
\begin{subequations}
\begin{equation}\label{def:Mp}
M_p(w,r)=\left( \meanint_{A^+_r} |w(x)|^p\,dx \right)^{1/p} \quad\hbox{where}\ A_r^+=\{x\in \RR^n_+: r<|x|<2r\}.
\end{equation}
Using this, we can also define
\begin{equation}\label{def:M1p}
M_{1,p}(w,r)=rM_p(\nabla w,r)+ M_p(w,r).
\end{equation}
\end{subequations}

\begin{proposition} \label{pr:potentialtheory} 
Suppose $F$ is the distribution \eqref{def:F} where $\vec f,f_0\in L^p_{\loc}(\overline{\RR^n_+}\backslash\{0\})$ for $p>n$ satisfy
\[
\int_{\{x\in\RR^n_+:|x|<1\}} \left(|\vec f(x)|+ |x f_0(x)|\right)|x|\,dx \, + \, 
\int_{\{x\in\RR^n_+:|x|>1\}} \left(|\vec f(x)|+|xf_0(x)|\right)\,|x|^{-1-n}\,dx\ <\ \infty.
\]
Then \eqref{eq:WeakSoln-Laplace-perp} defines a solution $w\in H^{1,p}_{\loc}(\overline{\RR^n_+}\backslash\{0\})$ of \eqref{eq:Neumann-Laplace-perp} that satisfies $Pw=0$ and
\[
M_{1,p}(w,\!r)\!\leq\! c\left(\! r^{-n}\!\!\int_0^r \!\left[ M_p(\vec f,\rho)\rho^{n}\!+\!M_p(f_0,\rho)\rho^{n+1}\right]\!d\rho
+ r^2\!\!\int_r^\infty\!\left[M_p(\vec f,\rho)\rho^{-2}\!+\!M_p(f_0,\rho)\rho^{-1}\right]\!d\rho \!\right).
\]
\end{proposition}

\noindent{\bf Proof.}
To obtain the desired estimates, let us assume $n\geq 3$, $r<|x|<2r$, and introduce the annulus $\widetilde A_r^+=\{x\in \RR^n_+: r/2<|x|<4r\}$. Then let us split the solution \eqref{eq:WeakSoln-Laplace-perp} into several parts:
\[
\begin{aligned}
w(x)&=\int_{\widetilde A^+_r}  \left(N(x,y) f_0(y)  -  \nabla_y N(x,y) \cdot  \vec f(y)\right) \, dy \\
-& \int_{r/2<|y|<|x|} \left( PN(x,y) f_0(y) - \nabla_y PN(x,y)\cdot\vec f(y) \right)\,dy\\
-& \int_{|x|<|y|<4r}  \left(PN(x,y) f_0(y) -  \nabla_y PN(x,y)\cdot\vec f(y) \right)\,dy \\
+&\int_{|y|<r/2}  \left( N^\perp(x,y) f_0(y)  - \nabla_y N^\perp(x,y) \cdot  \vec f(y) \right) \, dy \\
+& \int_{|y|>4r}  \left( N^\perp(x,y) f_0(y)  -  \nabla_y N^\perp(x,y) \cdot  \vec f(y)\right) \, dy\\
&=w_1(x)+w_2(x)+w_3(x)+ w_4(x)+w_5(x).
\end{aligned}
\]
(Here, and subsequently, by an integral such as $\int_{|y|<r/2}$ we actually mean the integral over $\{y\in\RR^n_+:|y|<r/2\}$.)
We estimate each of these terms separately.

The first term, $w_1$, can be estimated using classical results. For example, we can apply Theorem B* in \cite{SW} with $\lambda=n-1$, $\alpha=1$, $\beta=0$, and $p=q>n$ (which implies $p'<n$) to obtain
\[
\left\|\int_{\widetilde A_r^+} \nabla_y N(x,y)\cdot \vec f(y)\,dy\right\|_{L^p(A^+_r)}\leq c\,r\, \|\vec f\|_{L^p(\widetilde A^+_r)}.
\]
The same argument shows
\[
 r\left\|\int_{\widetilde A_r^+}  \frac{\partial}{\partial x_i}N(x,y)  f_0(y)\,dy\right\|_{L^p(A^+_r)}
\leq c\,r^2\, \|f_0\|_{L^p(\widetilde A^+_r)}.
\]
We  can also apply  Theorem B* in \cite{SW} with $\lambda=n-2$, $\alpha=2$, $\beta=0$, and 
$p=q>n$ to obtain
\[
\left\|\int_{\widetilde A_r^+}  N(x,y)  f_0(y)\,dy\right\|_{L^p(A^+_r)}
\leq c\,r^2\, \|f_0\|_{L^p(\widetilde A^+_r)}.
\]
Finally, we apply the $L^p$-boundedness of singular integral operators to obtain
\[
\left\|\int_{\widetilde A_r^+} \frac{\partial}{\partial x_i}\nabla_y N(x,y)\cdot \vec f(y)\,dy\right\|_{L^p(A^+_r)}\leq c\, \|\vec f\|_{L^p(\widetilde A^+_r)}.
\]
We conclude
\begin{equation}\label{est:w1}
M_{1,p}(w_1,r)\leq c\left(r\widetilde M_p(\vec f,r)+r^2\widetilde M_p(f_0,r)\right),
\end{equation}
where the tilde in $\widetilde M_p$ denotes that the spherical mean is taken over $\widetilde A^+_r$ instead of $A^+_r$.

For the second term, we note  $r/2<|y|<|x|<2r$ implies $|PN(x,y)|\leq c\,|x|^{2-n}\leq c\,r^{-n}|y|^2$ and $|\nabla_y PN(x,y)|\leq c\,|x|^{1-n}\leq c\, r^{-n}|y|$, so
\[
\left |\int_{r/2<|y|<|x|}  PN(x,y)  f_0(y)\,dy\right | \leq c\,r^{-n} \int_{r/2<|y|<|x|} |y|^2\,| f_0(y)|\,dy
\leq c\,r^{-n} \int_{|y|<2r} |y|^2 | f_0(y)|\,dy
\]
and
\[
\left |\int_{r/2<|y|<|x|} \nabla_y PN(x,y)\cdot \vec f(y)\,dy\right | \leq c\,r^{-n} \int_{r/2<|y|<|x|} |y||\vec f(y)|\,dy
\leq c\,r^{-n} \int_{|y|<2r} |y||\vec f(y)|\,dy.
\]
Similarly, we can estimate
\[
\left |\,r\int_{r/2<|y|<|x|} \frac{\partial}{\partial x_j} PN(x,y)  f_0(y)\,dy\right | \leq 
c\,r^{-n} \int_{|y|<2r} |y|^2 | f_0(y)|\,dy
\]
and
\[
\left |\,r\int_{r/2<|y|<|x|}\frac{\partial}{\partial x_j}  \nabla_y PN(x,y)\cdot \vec f(y)\,dy\right | 
\leq c\,r^{-n} \int_{|y|<2r} |y||\vec f(y)|\,dy.
\]
From these estimates we easily obtain
\begin{equation}\label{est:u2}
M_{1,p}(w_2,r)\leq c \,r^{-n} \int_{|y|<2r}\left( |y||\vec f(y)|+  |y|^2 | f_0(y|)\right)\,dy.
\end{equation}

For the third term we note $r<|x|<|y|<4r$ implies $|PN(x,y)|\leq c|y|^{2-n}\leq c\,r^{2}|y|^{-n}$ and $|\nabla_y PN(x,y)|\leq c|y|^{1-n}\leq c\,r^2|y|^{-n-1}$, so
\[
\left |\int_{|x|<|y|<4r}  PN(x,y)  f_0(y)\,dy\right | \leq c\,r^{2} \int_{|x|<|y|<4r} |y|^{-n}\,| f_0(y)|\,dy
\leq c\,r^{2} \int_{r<|y|} |y|^{-n} | f_0(y)|\,dy
\]
and
\[
\left| \int_{|x|<|y|<4r} \nabla_y PN(x,y)\cdot\vec f(y)dy\right| \leq c\,r^2  \int_{|x|<|y|<4r}|y|^{-n-1}\,|\vec f(y)|dy
\leq r^2 \int_{r<|y|}|y|^{-n-1}\,|\vec f(y)|dy.
\]
Similarly, we can estimate
\[
\left |\,r\int_{|x|<|y|<4r}  \frac{\partial}{\partial x_j} PN(x,y)  f_0(y)\,dy\right | 
\leq c\,r^{2} \int_{r<|y|} |y|^{-n} | f_0(y)|\,dy
\]
and
\[
\left| \,r\int_{|x|<|y|<4r} \frac{\partial}{\partial x_j} \nabla_y PN(x,y)\cdot\vec f(y)\,dy\right| 
\leq r^2 \int_{r<|y|}|y|^{-n-1}\,|\vec f(y)|\,dy.
\]
From these we easily obtain
\begin{equation}\label{est:u3}
M_{1,p}(w_3,r)\leq c \,r^{2} \int_{|y|>r}\left( |y|^{-n-1}|\vec f(y)|+  |y|^{-n} | f_0(y)|\right)\,dy.
\end{equation}

For the fourth term, we use $|y|<r/2<|x|/2<|x|$ to conclude $|N^\perp(x,y)|\leq c\,|y|^2/|x|^n\leq c\,r^{-n}|y|^2$
and $|\nabla_y N^\perp (x,y)|\leq c\,|y|/|x|^n\leq c\,r^{-n}|y|$. Consequently,
\[
\left|\int_{|y|<r/2} N^\perp (x,y) f_0(y) dy\right|\leq c \,r^{-n}\int_{|y|<r/2}|y|^{2}|f_0(y)|\,dy\leq
c \,r^{-n}\int_{|y|<r}|y|^{2}|f_0(y)|\,dy
\]
and
\[
\left|\int_{|y|<r/2} \nabla_yN^\perp (x,y) \cdot\vec f(y) dy\right|\leq c \,r^{-n}\int_{|y|<r/2}|y||\vec f(y)|\,dy\leq
c \,r^{-n}\int_{|y|<r}|y||\vec f(y)|\,dy.
\]
Similarly, we can estimate
\[
\left|\,r \int_{|y|<r/2} \frac{\partial}{\partial x_i}N^\perp (x,y) f_0(y) dy\right|\leq 
c \,r^{-n}\int_{|y|<r}|y|^2|f_0(y)|\,dy
\]
and
\[
\left|\,r\int_{|y|<r/2} \frac{\partial}{\partial x_i}\nabla_yN^\perp (x,y) \cdot\vec f(y) dy\right|\leq
c \,r^{-n}\int_{|y|<r}|y|\vec f(y)|\,dy.
\]
From these we easily obtain
\begin{equation}\label{est:u4}
M_{1,p}(w_4,r)\leq c \,r^{-n} \int_{|y|<r}\left( |y|^{2}|\vec f(y)|+  |y| | f_0(y)|\right)\,dy.
\end{equation}

For the fifth term, we use  $|x|<2r<4r<|y|$ to conclude $|N^\perp(x,y)|\leq c\,|x|^2/|y|^n\leq c\,r^2\,|y|^{-n}$ and $|\nabla_y N^\perp(x,y)|\leq c\,|x|^2/|y|^{n+1}\leq c\,r^2\,|y|^{-n-1}$. Hence
\[
\left| \int_{|y|>4r} N^\perp(x,y) f_0(y)\,dy\right|\leq c\,r^{2}\int_{|y|>4r} |y|^{-n}|f_0(y)|\,dy
\leq  c\,r^{2}\int_{|y|>r} |y|^{-n}|f_0(y)|\,dy
\]
and
\[
\left|\int_{|y|>4r} \nabla_yN^\perp (x,y) \cdot\vec f(y) dy\right|\leq c \, r^2 \int_{|y|>4r} |y|^{-n-1}|\vec f(y)|\,dy
\leq c \, r^2 \int_{|y|>r} |y|^{-n-1}|\vec f(y)|\,dy.
\]
Similarly, we estimate the first-order derivatives, so
we eventually obtain
\begin{equation}\label{est:u5}
M_{1,p}(w_5,r)\leq c\,r^2 \int_{|y|>r} \left(|\vec f(y)|+|y||f_0(y)|\right)|y|^{-n-1}\,dy.
\end{equation}

Putting these all together, we have
\[
\begin{aligned}
M_{1,p}(w,r)\leq c\left( r\widetilde M_p(\vec f,r)+r^2\widetilde M_p(f_0,r) 
+r^{-n}\int_{|y|<2r}\left(|\vec f(y)|+|y||f_0(y)|\right)\,|y|\,dy \right. \\
\left. +\,r^2\int_{|y|>r} \left(|\vec f(y)|+|y||f_0(y)|\right)|y|^{-n-1}\,dy\right).
\end{aligned}
\]
But
\[
r^{-n}\int_{r<|y|<2r}\left(|\vec f(y)|+|y||f_0(y)|\right)\,|y|\,dy \approx r^2\int_{r<|y|<2r} \left(|\vec f(y)|+|y||f_0(y)|\right)|y|^{-n-1}\,dy,
\]
so we can write this as
\begin{equation}\label{est-M1p-1}
\begin{aligned}
M_{1,p}(w,r)\leq c\left( r\widetilde M_p(\vec f,r)+r^2\widetilde M_p(f_0,r) 
+r^{-n}\int_{|y|<r}\left(|\vec f(y)|+|y||f_0(y)|\right)\,|y|\,dy \right. \\
\left. +\,r^2\int_{|y|>r} \left(|\vec f(y)|+|y||f_0(y)|\right)|y|^{-n-1}\,dy\right).
\end{aligned}
\end{equation}

Finally, the integrals in  \eqref{est-M1p-1} can be estimated in terms of $M_p$ and combined with the $\widetilde M_p$ term. For example, we can replace $|y|^2$ by $c\int_{|y|/2}^{|y|} \rho\,d\rho$, let $\rho=|z|$, and then interchange the order of integration to obtain
\[
\int_{|y|<r}\! |y|^2|f_0(y)|dy = c\!\int_{|y|<r}\int_{|y|/2<|z|<|y|} \!|z|^{2-n}|f_0(y)|dzdy
\leq c\!\int_{|z|<r}\!|z|^{2-n}\int_{|z|<|y|<2|z|}\!|g(y)|dydz.
\]
But, by the H\"older inequality,
\[
\begin{aligned}
\int_{|z|<|y|<2|z|}|f_0(y)|dy & \leq \left( \int_{|z|<|y|<2|z|}|f_0(y)|^pdy\right)^{1/p}\left(\int_{|z|<|y|<2|z|}\,dy\right)^{1/p'} \cr
& =c  \left( \int_{|z|<|y|<2|z|}|g(y)|^pdy\right)^{1/p} |z|^{n/p'}=c\,M_p(f_0,|z|)\,|z|^n.
\end{aligned}
\]
Thus
\[
\int_{|y|<r}\! |y|^2|f_0(y)|\,dy \leq c\int_{|z|<r} |z|^2\,M_p(f_0,|z|)\,dz=c\int_0^r \rho^{n+1}\,M_p(f_0,\rho)\,d\rho.
\]
Similarly, we can show
\[
\int_{|y|>r} |y|^{-n} |f_0(y)| \,dy \leq c\,\int_r^\infty \rho^{-1} M_p(f_0,\rho)\,d\rho.
\]
If we similarly estimate the analogous integrals involving $\vec f$, we will obtain the estimate in the proposition.
 $\Box$

\section{Variable Coefficients in the Half-space Problem}

In this section we consider (\ref{eq:weakform}) when $U=\RR^n_+$, i.e.\ we assume 
that  $u\in H^{1,2}_\loc(\overline{\RR^n_+})$  satisfies
\begin{equation}\label{eq:HalfspaceProblem}
\int_{\RR^n_+}  a_{ij}\,\partial_j u\, \partial_i\eta \,dx=0
\quad\hbox{for all}\ \eta\in C_{comp}^1(\overline{\RR^n_+}).
\end{equation}
We want to consider the regularity of $u$ at a point on $\RR^{n-1}=\partial\RR^{n}_+$ which, for convenience, we take to be the origin. As shown in Appendix \ref{SobolevReg}, the continuity of the $a_{ij}$ enables us to conclude that $u\in H^{1,p}_{\loc}(\overline{\RR^n_+})$ for all $p>2$. Let us fix $p\in (n,\infty)$. 
By a change of independent variables we may arrange $a_{ij}(0)=\delta_{ij}$, so we assume that the coefficients satisfy 
\begin{equation}\label{aij-deltaij}
\sup_{|x|=r} |a_{ij}(x)-\delta_{ij}| \leq \om(r)\quad\hbox{as}\ r\to 0,
 \end{equation}
 where $\om$ is a continuous, nondecreasing function satisfying \eqref{eq:Sq-Dini} and 
 \eqref{om-vanishing}. 
We shall also assume that we have scaled the independent variables so that for $\delta$ very small we have
\begin{equation}\label{om-delta}
\int_0^1\frac{\om^2(r) }{r}\,dr <\delta
\ \hbox{and}\ 
\om(1)=\delta
\end{equation}
For convenience, we extend $\om$ to satisfy $\om(r)=\delta$ for $r>1$.

Now let us introduce a smooth cut-off function $\chi(r)$ satisfying $\chi(r)=1$ for $r<1/4$ and $\chi(r)=0$ for $r>3/4$.
Then $\chi(|x|) u(x)$ is a compactly supported function that agrees with $u(x)$ near $x=0$. What equation does $\chi u$ satisfy?
If we replace $\eta$ in \eqref{eq:weakform} by $\chi\eta$ and rearrange, we obtain
\[
\int_{\RR^n_+} \left( a_{ij}\,\partial_j(\chi u)\, \partial_i\eta 
- f_i\partial_i\eta+  f_0\eta\right)\,dx=0
\quad\hbox{for all}\ \eta\in C_{comp}^1(\overline{\RR^n_+}),
\]
where
$ f_i:=a_{ij} u\, \partial_j\chi$ and $ f_0:=a_{ij}\partial_j u\, \partial_i\chi$
are known to be in $L^p_{\rm comp}(\overline{\RR^n_+})$.
Since we are interested in the behavior of $u$ near $x=0$ where $u$ and $\chi u$ agree, after relabeling we can assume that $u\in H^{1,p}(\RR^n_+)$ has support in $|x|<1$ and satisfies
\begin{equation}\label{eq:HalfspaceProblem1}
\int_{\RR^n_+} \left( a_{ij}\,\partial_j u\, \partial_i\eta 
-f_i\partial_i\eta+ f_0\eta\right)\,dx=0
\quad\hbox{for all}\ \eta\in C_{comp}^1(\overline{\RR^n_+}),
\end{equation}
where $f_i,f_0\in L^p$ have support in $|x|<1$ with $f_i(r)=0$ for $r<1/4$. Since $u$ vanishes outside $|x|<1$, there is no harm in assuming that $a_{ij}$ satisfies
\begin{equation}
a_{ij}(x)=\de_{ij} \quad\hbox{ for}\ |x|\geq 1. 
\label{eq:aij_for_r>2}
\end{equation}

Let us recall the  decomposition $u(x)=u_0(r)+\widetilde v(r)\cdot\widetilde x+w(x)$ as defined in \eqref{spectral-decomp}. Since we  assumed that $u$ is supported in $|x|<1$, we have that $u_0,\widetilde v$, and $w$ are all supported in  $|x|<1$. Moreover, as shown in Appendix B,
\begin{equation}\label{grad(u)-props}
\nabla u\in L^2(B_+(1)) \ \Rightarrow\ 
\int_0^1 \left[ (u'_0)^2+|\widetilde v|^2 +r^2|\widetilde v\,'|^2\right] r^{n-1}\,dr <\infty 
\ \ \hbox{and}\ \ \nabla w\in L^2(B_+(1)).
\end{equation}
 
 To formulate the connection between the decomposition and the dynamical system,  let  $r=e^{-t}$ and introduce
 \begin{subequations}
\begin{equation}\label{def:epsilon}
\e(t):=\om(e^{-t}) \quad\hbox{for $-\infty<t<\infty$.}
\end{equation}
Notice that
\begin{equation}
\int_0^\infty \e^2(t)\,dt=\int_0^1\frac{\om^2(r)}{r}\,dr.
\end{equation}
\end{subequations}
To control the behavior of $\widetilde v(r)$ and $r\widetilde v\,'(r)$ as $r\to 0$, we need to control the behavior of $\widetilde v(t)$ and $\widetilde v_t(t)$ as $t\to \infty$.
In Appendix D of this paper, we show that new dependent variables $(\varphi,\psi)$ can be introduced that satisfy a $2(n-1)$-dimensional dynamical system
 \begin{subequations}\label{ourODEsystem:joint}
 \begin{equation} \label{ourODEsystem:a}
\frac{d}{dt}
\begin{pmatrix}
\varphi \\ \psi
\end{pmatrix}
+
\begin{pmatrix}
0 & 0 \\ 0 & -nI
\end{pmatrix}
\begin{pmatrix}
\varphi \\ \psi
\end{pmatrix}
+
{\mathcal R}(t)
\begin{pmatrix}
\varphi \\ \psi
\end{pmatrix}
=  
g(t,\nabla w)+h(t) \ \hbox{for}\ T<t<\infty,
\end{equation}
where the matrix ${\mathcal R}$ depends upon the coefficients $a_{ij}$ and can be decomposed into blocks
\begin{equation}\label{ourODEsystem:b}
{\mathcal R}(t)=
\begin{pmatrix}
R_{1}(t) & R_{2}(t) \\ R_{3}(t) & R_{4}(t)
\end{pmatrix}
\quad
\hbox{with}\ \pmb{\bigr |}R_j(t)\pmb{\bigr |}\leq \e(t).
\end{equation}
The block $R_1$ satisfies 
\begin{equation}\label{ourODEsystem:c}
|R_1(t)-R(t)|\leq c\,\e^2(t) \ \hbox{as}\ t\to\infty, 
\end{equation}
where the 
$(n-1)\times(n-1)$ matrix $R(t)$ is given by \eqref{def:R-halfspace}.
The term  $g(t,\nabla w)$ in \eqref{ourODEsystem:a} denotes a vector function of $t$ that depends on $\nabla w$ (the gradient in the $x$-variables) in such a way that 
\begin{equation}
|g(t,\nabla w)|\leq c\, \e(t)\,\meanint_{S_+^{n-1}}  |\nabla w|\,ds\quad\hbox{for}\ t>0,
 \label{ourODEsystem:d}
\end{equation}
and the term $h$ in \eqref{ourODEsystem:a} is a vector function in $L^1(0,\infty)$ with $L^1$-norm satisfying
\begin{equation} \label{ourODEsystem:e}
\|h\|_1\leq  c\,\left(\|\vec f\|_p+\|f_0\|_p\right).
\end{equation}
 Moreover, the difference between the new dependent variables $(\varphi,\psi)$ and $(\widetilde v,\widetilde v_t)$ is estimated by
 \begin{equation} \label{ourODEsystem:f}
\left| \begin{pmatrix}
\widetilde v \\
\widetilde v_t
\end{pmatrix}
-
 \begin{pmatrix}
n (\varphi + \psi) \\ 
n^2\psi
\end{pmatrix}
\right|
\leq
c\,\e(t)\left(
|\varphi(t)|+|\psi(t)|+\meanint|\nabla w|ds
\right).
\end{equation}
\end{subequations}
We will use this and the stability of $(\varphi,\psi)$ as $t\to\infty$ to control the behavior of $\widetilde v$ as $r\to 0$.

With these preliminaries, we are able to prove the following.
\begin{theorem} \label{th:1} 
Suppose the $a_{ij}$ satisfy \eqref{aij-deltaij} where $\om$ satisfies \eqref{eq:Sq-Dini} and 
 \eqref{om-vanishing}, and  the dynamical system \eqref{eq:DynSystm} with matrix $R$ given by \eqref{def:R-halfspace} is uniformly stable as $t\to\infty$. Then every weak solution $u\in H_\loc^{1,2}(\RR^n_+)$ of \eqref{eq:HalfspaceProblem} is Lipschitz continuous at $x=0$. If, in addition, every solution of the dynamical system 
\eqref{eq:DynSystm} is asymptotically constant as $t\to\infty$, then $u$ is differentiable at $x=0$, and
\[
\partial_j u(0)=
\lim_{r\to 0}\frac{n}{r}\meanint_{S^{n-1}_+} u(r\theta)\,\theta_j\,ds \quad \hbox{for}\ j=1,\dots,n-1, \qquad
\partial_n u(0)= 0.
\]
\end{theorem}

\noindent{\bf Proof.} As indicated above, we may assume for some $p\in (n,\infty)$ that $u\in H^{1,p}(\overline{\RR^n_+})$ is supported in $|x|<1$ and satisfies \eqref{eq:HalfspaceProblem1}. The strategy of the proof is to construct a solution $u^*$ of \eqref{eq:HalfspaceProblem1} in the
form \eqref{spectral-decomp}.  This is done by finding $w$ as a fixed point for a certain map $S$ on the Banach space $Y$, which is defined to be $w\in H^{1,p}_{\loc}(\overline{\RR^n_+}\backslash\{0\})$ with finite norm
\begin{equation}
\| w \|_Y = \sup_{0<r<1} \frac{M_{1,p}(w,r)}{\om(r)\, r} + \sup_{r>1}\frac{M_{1,p}(w,r)}{r^{-n}}.
\end{equation}
Since $p>2$ we see that $w\in Y$ implies $M_2(\nabla w,r)\leq C\,\om(r)$ for $0<r<1$. 
 As we shall see, finding $w$ also yields $\widetilde v$ and $r\widetilde v_r$ from the solution of the dynamical system \eqref{ourODEsystem:a}.
Moreover, $u_0'$ can be found in terms of $\widetilde v$, $r\widetilde v_r$, and $w$, and we find that the stability properties of the dynamical system \eqref{ourODEsystem:a} control the asymptotic behavior of $\widetilde v$ and $r\widetilde v_r$ as $r\to 0$, and hence also of $u_0$. Under the assumed stability of \eqref{ourODEsystem:a}, the constructed $u^*$ has the required regularity, and it only remains to show that $u^*=u$; this is done using the uniqueness of solutions of \eqref{eq:HalfspaceProblem1} discussed in Appendix E. Let us now discuss the details of this argument.

For a given $w\in Y$, we want to solve \eqref{ourODEsystem:a} with initial conditions $\phi(0)=0=\psi(0)$ to find $(\phi,\psi)$ and hence $\widetilde v$, $r \widetilde v_r$. To control the dependence of $\widetilde v$ on $w$, let us write $\widetilde v=\widetilde v^w+\widetilde v^0$ where $\widetilde v^w$ corresponds to solving \eqref{ourODEsystem:a} with $h\equiv 0$ and $\widetilde v^0$ corresponds to solving it with $g(t,\nabla w)\equiv 0$. In order to estimate $\widetilde v^w$ on $(0,\infty)$, we will use Proposition \ref{pr:stabilty} in Appendix E. Consequently, we need $g=(g_1,g_2)$ to satisfy: i)  $g_1\in L^1(0,\infty)$ and ii) $g_2$ satisfies \eqref{g2-condition}. First, we use \eqref{ourODEsystem:d} to conclude
\begin{subequations}\label{est:g1,alpha}
\begin{equation}\label{est:g1}
\begin{aligned}
\int_0^\infty |g_1(t,\nabla w)|\,dt  & \leq c\,\left(\int_0^\infty \e^2(t)\,dt\right)^{1/2}\left(\int_0^\infty\int_{S_+^{n-1}}|\nabla w|^2\,ds\,dt\right)^{1/2}\\
& \leq c\,\sqrt{\delta}\left(\int_0^1\meanint_{S_+^{n-1}}|\nabla w|^2 ds\,\frac{d\rho}{\rho}\right)^{1/2}.
\end{aligned}
\end{equation}
We will conclude the finiteness of this bound below.
Second, we use \eqref{ourODEsystem:d} to conclude
\[
e^{\alpha t} \int_t^\infty |g_2(\tau,\nabla w)|\,e^{-\alpha\tau}\,d\tau 
\leq c\,\e(t)\int_t^\infty e^{\alpha(t-\tau)}\meanint_{S^{n-1}_+}|\nabla w|\,ds
\leq c_\alpha\,\e(t), 
\]
where
\begin{equation}
c_\alpha=\frac{c}{\sqrt{2\alpha}}\left(\int_0^1\meanint |\nabla w|^2ds\frac{d\rho}{\rho}\right)^{1/2}.
\end{equation}
Now let us perform a calculation for $1\leq p<\infty$: for $j=0,1,\dots$, let $r_j=2^{-j}$ so
\begin{equation}\label{int-nabla_w-est}
\begin{aligned}
\int_0^1\int_{S^{n-1}_+}|\nabla w|^p\,ds\frac{d\rho}{\rho}=\sum_{j=1}^\infty \int_{r_j}^{2r_j}\int_{S^{n-1}_+} |\nabla w|^p\,ds\,\frac{d\rho}{\rho}
\leq c\sum_{j=1}^\infty M_p^p(\nabla w,r_j).
\end{aligned}
\end{equation}
As observed above, $M_2(\nabla w,r)\leq C\,\om(r)$ as $r\to 0$, so we may apply
the above estimate with $p=2$ and the following calculation
\begin{equation}
\sum_{j=1}^\infty \om^2(r_j)=2\sum_{j=1}^\infty \om^2(r_j)\frac{r_{j-1}-r_{j}}{r_{j-1}}\leq 2\int_0^1\om^2(\rho)\frac{d\rho}{\rho}< 2\delta
\end{equation}
\end{subequations}
to conclude the finiteness of $c_\alpha$ and the bound \eqref{est:g1}. Thus we have confirmed i) and ii).

Next, let us describe the variational PDE that $w$ satisfies. 
As in \cite{MM3} we introduce $\Omega_{ij}=a_{ij}-\delta_{ij}$, which
satisfies $|\Omega_{ij}(x)|\leq \om(r)$ for $0<r=|x|<1$ and $\Omega_{ij}(x)=0$ for $|x|>1$. 
Now the variational problem (\ref{eq:HalfspaceProblem1}) can be written as
\[
\int_{\RR^n_+}\left( \nabla u\cdot\nabla \eta +\Omega_{ij}\partial_j u\, \partial_i \eta-f_i\partial_i\eta +f_0\eta \right)dx=0
\quad\hbox{for all}\ \eta\in C_{comp}^1(\overline{\RR^n_+}).
\]
Since this holds for all $\eta$, it holds for $\eta^\perp$:
\[
\int_{\RR^n_+} \left(\nabla u\cdot\nabla (\eta^\perp) +\Omega_{ij}\partial_j u\, \partial_i (\eta^\perp)-f_i\partial_i(\eta^\perp) +f_0\,\eta^\perp \right)dx=0.
\]

Now we claim that
\begin{equation}\label{eq:int-grad(Pu)grad(eta-perp)}
\int_{\RR^n_+} \nabla (Pu)\cdot\nabla (\eta^\perp)\,dx=0= \int_{\RR^n_+} \nabla w\cdot\nabla (P\eta) \,dx
\quad\hbox{for all}\ \eta\in C_{comp}^1(\overline{\RR^n_+}).
\end{equation}
To prove this, let us first assume $u\in C^2_{\rm comp}(\overline{\RR^n_+})$. Then, by Lemma \ref{le:d_nPu=0}, we have  $Pu\in C^2_{\rm comp}(\overline{\RR^n_+})$, $P\eta\in C^1_{\rm comp}(\overline{\RR^n_+})$, and $\partial Pu/\partial x_n=0=\partial P\eta/\partial x_n$ on $\RR^{n-1}$. Applying the divergence theorem, we conclude
\[
\int_{\RR^n_+} \nabla (Pu)\cdot\nabla (\eta^\perp)\,dx=- \int_{\RR^n_+}\Lap (Pu) \,\eta^\perp\,dx=- \int_{\RR^n_+}(\Lap (Pu))^\perp \,\eta\,dx.
\]
However, for fixed $r$, $Pu(r,\cdot)\in V=$ span$(1,\theta_1,\dots,\theta_{n-1})$ on $S^{n-1}_+$, and 
 $\Lap$ preserves $V$, so $(\Lap (Pu))^\perp=0$; this proves the first equality in \eqref{eq:int-grad(Pu)grad(eta-perp)}
 when  $u\in C^2_{\rm comp}(\overline{\RR^n_+})$. In general, for $u\in H^{1,p}_{comp}(\overline{\RR^n_+})$,
 we extend $u$  by zero to $\RR^n_-$, and mollify by $u_\e=\phi_\e\star u$, where
 $\phi_\e(x)=\e^{-n}\phi(|x|/\e)$ with $\phi\in C^\infty(\RR^n)$ satisfying supp($\phi$)$\subset B_1(0)$ and $\int \phi(x)\,dx=1$. Then  $u_\e$  is smooth (on all of $\RR^n$), and we can apply the above argument to conclude 
\[
\int_{\RR^n_+} \nabla (P\,u_\e)\cdot\nabla (\eta^\perp)\,dx=0.
\]
Since we have assumed 
$u\in H^{1,p}(\overline{\RR^n_+})$, we can show in the standard way that $u_\e\to u$ in $ H^{1,p}(\overline{\RR^n_+})$ as $\e\to 0$. Then we can use Lemma \ref{le:Puj->Pu} to conclude (even without the $\perp$ on $\eta$)
that
\[
\int_{\RR^n_+}\nabla(P\,u_\e)\cdot\nabla\eta\,dx
\to \int_{\RR^n_+}\nabla(Pu)\cdot\nabla\eta\,dx \quad\hbox{as}\ \e\to 0 \quad\hbox{for all $\eta\in C^1_{comp}(\overline{\RR^n_+}\backslash\{0\})$}.
\]
This establishes the first equality in  \eqref{eq:int-grad(Pu)grad(eta-perp)}.
The second equality in  \eqref{eq:int-grad(Pu)grad(eta-perp)} follows a similar argument. 

But \eqref{eq:int-grad(Pu)grad(eta-perp)}  means that
\begin{equation}\label{eq:int-grad(u)grad(eta-perp)}
\int_{\RR^n_+} \nabla u\cdot\nabla( \eta^\perp)\,dx
=\int_{\RR^n_+} \nabla w\cdot\nabla \eta\,dx
\quad\hbox{for all}\ \eta\in C_{comp}^1(\overline{\RR^n_+}).
\end{equation}
Using this and the fact that $u_0'$ can be expressed in terms of $\widetilde v^w$ and $w$ (see \eqref{eqnfor-u0'} in Appendix D), we see that the variational problem that $w$ satisfies can be written as
\begin{equation}\label{generalvariational3}
\int_{\RR^n_+} \nabla w\cdot\nabla \eta \,dx+F_{1,w}^\perp[\eta]+F_{1,0}^\perp[\eta]+F_0^\perp[\eta]=0 \quad\hbox{for all $\eta\in C_{comp}^{1}(\overline{\RR^n_+})$},
\end{equation}
where
\begin{equation}
F_{1,w}^\perp[\eta]=\int_{\RR^n_+} \vec f^{\,w}\cdot\nabla(\eta^\perp) dx, \quad
F_{1,0}^\perp[\eta]=\int_{\RR^n_+} \vec f^{\,0}\cdot\nabla(\eta^\perp) dx, \quad
F_0^\perp[\eta]=\int_{\RR^n_+} f_0\eta^\perp dx,
\end{equation}
with the vector functions $\vec f^{\,w}$ and $\vec f^{\,0}$ defined by
\begin{equation}
f_i^{w}=\Omega_{ij}\left(\partial_j w - 
\frac{ r\widetilde\beta\cdot(\widetilde v^w)'+\widetilde\gamma\cdot\widetilde v^w+p[\nabla w]}{\alpha}\theta_j
+\partial_j(\widetilde x\cdot\widetilde v^w)\right)
\end{equation}
\begin{equation}
f_i^{0}=\Om_{ij}\left(\frac{ \vartheta(r)-r\widetilde\beta\cdot(\widetilde v^0)'-\widetilde\gamma\cdot\widetilde v^0}{\alpha}\theta_j
+\partial_j(\widetilde x\cdot\widetilde v^0)\right)
-f_i.
\end{equation}
Here, as in Appendix D, the functions $\alpha$, $\tilde\beta$, $\tilde\gamma$, $p[\nabla w]$, and $\vartheta$ of $r$ satisfy
\begin{equation}\label{coefficients,r<1}
\begin{aligned}
|\alpha(r)-1|, |\widetilde\beta(r)|, |\widetilde\gamma(r)| \leq\om(r) \quad \hbox{for}\ 0<r<1,
\\
|p[\nabla w](r)|\leq \om(r) \meanint_{S_+^{n-1}} |\nabla w(r\theta)|\,ds  \quad \hbox{for}\ 0<r<1,
\\
|\vartheta(r)|\leq \meanint_{S_+^{n-1}} (|\vec f(r\theta)|+|f_0(r\theta)|)\,ds \quad \hbox{for}\ 0<r<1,
\end{aligned}
\end{equation}
and
$\alpha(r)=1$ and $\widetilde\beta(r)=\widetilde\gamma(r)=p[\nabla w](r)=\vartheta(r)=0$ for $r>1$.

For $w\in Y$, define $z=S(w)$ to be the solution of 
\begin{equation}\label{eq:variational-z}
\int_{\RR^n_+} \nabla z\cdot\nabla \eta \,dx+F_{1,w}^\perp[\eta]+F_{1,0}^\perp[\eta]+ F^\perp_{0}[\eta]=0
\quad\hbox{for all $\eta\in C_{comp}^{1}(\overline{\RR^n_+})$}
\end{equation}
\noindent
that is provided by Proposition \ref{pr:potentialtheory}, i.e.\
\begin{equation}
S(w)=z(x)=\int_{\RR^n_+}\left(N^\perp(x,y)f_0(y)-\nabla_y N^\perp(x,y)\cdot \vec f^0(y)
-\nabla_y N^\perp(x,y)\cdot \vec f^w(y)
\right)dy.
\end{equation}
If we can show that $S:Y\to Y$ has a fixed point $w$, then this is the solution of \eqref{generalvariational3} that we seek.

To show $S$ has a fixed point, we write 
$
Sw=\xi-Tw
$
where
\begin{equation}\label{def:xi}
\xi(x)= \int_{\RR^n_+}\left(N^\perp(x,y)f_0(y)-\nabla_y N^\perp(x,y)\cdot \vec f^0(y)\right)dy
\end{equation}
and
\begin{equation}\label{def:T}
Tw(x)=\int_{\RR^n_+}\nabla_y N^\perp(x,y)\cdot \vec f^w(y)
\, dy.
\end{equation}
If we can show that $\xi\in Y$ and $T:Y\to Y$ with small norm, then we can solve $w+Tw=\xi$ to find our fixed point $w=Sw$. 
To estimate $M_{1,p}(Tw,r)$ we will apply Proposition \ref{pr:potentialtheory}:
\begin{equation}\label{est:M(Tw)}
M_{1,p}(Tw,\!r)\!\leq\! c\left(\! r^{-n}\!\!\int_0^r   M_p(\vec f^w,\rho)\rho^{n} d\rho
+ r^2\!\!\int_r^\infty M_p(\vec f^w,\rho)\rho^{-2}\, d\rho \!\right).
\end{equation}
So we only need to estimate $M_p(\vec f^w,r)$ and integrate.

Now let us assume $\|w\|_Y\leq 1$ and show that $\|Tw\|_Y$ is small. We split $Tw$ into three terms:
\[
T_1w(x)=\int_{\RR^n_+}\nabla_y N^\perp(x,y)\cdot \Om\nabla w(y) \, dy
\]
\[
T_2w(x)=\int_{\RR^n_+}\nabla_y N^\perp(x,y)\cdot \Om\nabla(\tilde y\cdot \tilde v^w) \, dy
\]
\[
T_3w(x)=\int_{\RR^n_+}\frac{1}{\alpha(r_y)}\left(r_y\,\widetilde\beta\cdot(\widetilde v^w)'(r_y)+\widetilde\gamma\cdot\widetilde v^w(r_y)+p[\nabla w](r_y)\right)\nabla_y N^\perp(x,y)\cdot\Om\theta(y)\,dy.
\]
Here $r_y:=|y|$ and we have written the vector $\Om_{ij}\partial_j w$ simply as $\Om\nabla w$; similarly for
 $ \Om\nabla(\tilde y\cdot \tilde v^w)$ and $\Om\theta$.
Let us first consider $T_1w$. Recall that $|\Om(r)|\leq \om(r)$ for $0<r<1$ and $\Om(r)\equiv 0$ for $r>1$, so
\[
\begin{aligned}
r^{-n}\!\!\int_0^r   & M_p(\Om\nabla w,\rho)\rho^{n} d\rho
 + r^2\!\!\int_r^\infty M_p(\Om\nabla w,\rho)\rho^{-2}\, d\rho \\
& \leq \begin{cases} r^{-n}\!\!\int_0^r   \om(\rho) M_p(\nabla w,\rho)\rho^{n} d\rho
+ r^2\!\!\int_r^1 \om(\rho)M_p(\nabla w,\rho)\rho^{-2}\, d\rho & \hbox{for $0<r<1$} \\
r^{-n}\int_0^1 \om(\rho) M_p(\nabla w,\rho)\rho^n\,d\rho & \hbox{for $r>1$.}
\end{cases}
\end{aligned}
\]
For $0<r<1$ we have assumed $\|w\|_Y\leq 1$ and we have $M_{1,p}(w,r)\leq \om(r)\,r$, so $M_p(\nabla w,r)\leq \om(r)$ and we can estimate
\[
\begin{aligned}
r^{-n}\!\int_0^r &  \om(\rho) M_p(\nabla w,\rho)\rho^{n} d\rho + r^2\!\int_r^1 \om(\rho)M_p(\nabla w,\rho)\rho^{-2}\, d\rho \\
&\leq r^{-n}\int_0^r \om^2(\rho)\rho^n\,d\rho+r^2\delta\int_r^1\om(\rho)\rho^{-2}\,d\rho\leq c\left( \om^2(r)r+ \delta\om(r)r(r^\kappa-1)\right)\leq c\,\delta\,\om(r)\,r.
\end{aligned}
\]
Consequently, for $0<r<1$ we have from \eqref{est:M(Tw)} that
\[
M_{1,p}(T_1w,\!r) \leq c\,\delta\,\om(r)\,r.
\]
Meanwhile, for $r>1$, \eqref{est:M(Tw)} implies
\[
\begin{aligned}
M_{1,p}(T_1w,r) & \leq r^{-n}\int_0^1 \om(\rho) M_p(\nabla w,\rho)\rho^n\,d\rho \\
& \leq r^{-n}\int_0^1 \om^2(\rho)\rho^n\,d\rho \leq c\,\delta\, r^{-n}.
\end{aligned}
\]
Thus $\|T_1w\|_Y\leq c\,\delta$ and, if we take $\delta$ sufficiently small, we can arrange that $T_1:Y\to Y$ has norm less than $1/3$.

Next consider $T_2 w$. Again we use $|\Om(r)|\leq \om(r)$ for $0<r<1$ and $\Om(r)\equiv 0$ for $r>1$ to obtain
\[
\begin{aligned}
r^{-n}\!\!\int_0^r   & M_p(\Om\nabla (\tilde y\cdot \tilde v^w),\rho)\rho^{n} d\rho
 + r^2\!\!\int_r^\infty M_p(\Om\nabla (\tilde y\cdot \tilde v^w),\rho)\rho^{-2}\, d\rho \\
& \leq \begin{cases} r^{-n}\!\!\int_0^r   \om(\rho) M_p(\nabla (\tilde y\cdot \tilde v^w),\rho)\rho^{n} d\rho
+ r^2\!\!\int_r^1 \om(\rho)M_p(\nabla (\tilde y\cdot \tilde v^w),\rho)\rho^{-2}\, d\rho & \hbox{for $0<r<1$} \\
r^{-n}\int_0^1 \om(\rho) M_p(\nabla (\tilde y\cdot \tilde v^w),\rho)\rho^n\,d\rho & \hbox{for $r>1$.}
\end{cases}
\end{aligned}
\]
To estimate $\nabla (\tilde y\cdot \tilde v^w)$ we need to estimate $\tilde v^w$ and $r(\tilde v^w)'$. But, using \eqref{ourODEsystem:f}, these can be expressed in terms of the solution $(\phi,\psi)$ of the dynamical system \eqref{ourODEsystem:a}.
Thus we find
\[
\sup_{|y|<1}|\nabla(\tilde y\cdot \tilde v^w)|\leq c\sup_{r<1} (r|(\tilde v^w)'|+|\tilde v^w|)
\leq c\sup_{t>0} (|\phi(t)|+|\psi(t)|)\leq c(c_\alpha+\|g_1\|_1),
\]
where we have used Proposition \ref{pr:stabilty}  in Appendix E for the last estimate. Now we can estimate $c_\alpha$ and $\|g_1\|_1$ as in \eqref{est:g1,alpha} to find $c_\alpha \leq c\,\sqrt{\delta}$ and $\|g_1\|_1\leq c\,\delta$. So we conclude that for $0<\rho<1$
\[
M_p(\nabla (\tilde y\cdot \tilde v^w),\rho) \leq c \sup_{|y|<1}|\nabla(\tilde y\cdot \tilde v^w)|\leq c\,\sqrt{\delta}.
\]
We can use this in \eqref{est:M(Tw)} to estimate $M_{1,p}(T_2 w,r)$:
\[
M_{1,p}(T_2 w,r)\leq 
\begin{cases}
c\left( r^{-n}\int_0^r \om(\rho) \sqrt{\delta}\rho^n\,d\rho+r^2\int_r^1\om(\rho)\sqrt{\delta}\rho^{-2}\,d\rho\right)\leq c\,\sqrt{\delta}\,\om(r)\,r
& \hbox{for $0<r<1$}\\
c\,r^{-n}\int_0^1 \om(\rho)\sqrt{\delta}\rho^n\,d\rho\leq c\,\delta^{3/2}\,r^{-n}  & \hbox{for $r>1$.}
\end{cases}
\]
Thus, $\|T_2 w\|_Y\leq c\,(\sqrt{\delta}+\delta^{3/2})$ and, if we take $\delta$ sufficiently small, we can arrange that $T_2:Y\to Y$ has norm less than $1/3$.

Finally we consider $T_3w$. We first need to estimate $M_p(\alpha^{-1}(r\widetilde\beta\cdot(\widetilde v^w)'+\widetilde\gamma\cdot\widetilde v^w+p[\nabla w])\Om\theta,r)$  for $0<r<1$. But, recalling the properties \eqref{coefficients,r<1}
and some of the estimates used for $T_2$, we have
\[
\begin{aligned}
M_p(\alpha^{-1}(r\widetilde\beta\cdot(\widetilde v^w)'+\widetilde\gamma\cdot\widetilde v^w+p[\nabla w])\Om\theta,r)
&\leq c\,\om^2(r)\left[M_p(r(\widetilde v^w)',r)+M_p(\widetilde v^w,r)+M(\nabla w,r)\right] \\
&\leq c\,\sqrt{\delta}\,\om^2(r).
\end{aligned}
\]
Applying Proposition \ref{pr:potentialtheory}, we obtain for $0<r<1$
\[
M_{1,p}(T_3w,r)\leq c\left(r^{-n}\int_0^r \sqrt{\delta}\,\om^2(\rho)\rho^n\,d\rho + r^2\int_r^1\sqrt{\delta}\,\om^2(\rho)\rho^{-2}\,d\rho\right)
\leq c\,\delta^{3/2}\om(r)\,r.
\]
Meanwhile, for $r>1$ we simply have
\[
M_{1,p}(T_3w,r)\leq c\,r^{-n}\int_0^1\sqrt{\delta}\,\om(\rho)\rho^n\,d\rho
\leq c\,\delta^{3/2}\, r^{-n}.
\]
Thus $\|T_3w\|_Y\leq c\,\delta^{3/2}$, and if we take $\delta$ sufficiently small, we can arrange $T_3:Y\to Y$ to have norm less than $1/3$. Consequently, $T=T_1+T_2-T_3:Y\to Y$ has norm less than $1$.

To show that $\xi$ defined in \eqref{def:xi} is in $Y$, let us split it up into several terms: $\xi=\xi_1-\xi_2+\xi_3+\xi_4$, where
\[
\begin{aligned}
\xi_1(x)&=\int_{\RR^n_+} N^\perp(x,y)\,f_0(y)\,dy, \\
\xi_2(x)&=\int_{\RR^n_+} \nabla_y N^\perp(x,y)\cdot\vec f(y)\,dy, \\
\xi_3(x)&=\int_{\RR^n_+} \nabla_y N^\perp(x,y)\cdot\Om\nabla(\widetilde x\cdot\widetilde v^0)\,dy, \\
\xi_4(x)&=\int_{\RR^n_+} \frac{1}{\alpha(r_y)}\left(\vartheta(r_y)-r_y\,\widetilde\beta\cdot(\widetilde v^0)'(r_y)
-\widetilde\gamma\cdot \widetilde v^0(r_y) \right)\nabla_yN^\perp(x,y)\cdot\Om\theta\,dy.
\end{aligned}
\]
Since $f_0\in L^p$ and is supported in $|x|<1$, we can apply Proposition \ref{pr:potentialtheory}:
\[
\begin{aligned}
M_{1,p}(\xi_1,r) & \leq c\left(r^{-n}\int_0^r M_p(f_0,\rho)\rho^{n+1}\,d\rho+r^2\int_r^1 M_p(f_0,\rho)\,d\rho\right)\\
& \leq \begin{cases}
c\,r^2\,\| f_0\|_p  & \hbox{for}\ 0<r<1 \\
c\,r^{-n}\,\| f_0\|_p  & \hbox{for}\ r>1.
\end{cases}
\end{aligned}
\]
Since $r\leq\om(r)$ for $0<r<1$, we see that $\xi_1\in Y$. To estimate $\xi_2$ recall that $\vec f$ is supported in $1/4<|x|<1$, so 
$|\vec f(x)|\leq c\,\om(|x|)\,|\vec f(x)|$, and we obtain from Proposition \ref{pr:potentialtheory} the estimate
\[
\begin{aligned}
M_{1,p}(\xi_2,r) & \leq c\left(r^{-n}\int_0^r \om(\rho)\,M_p(\vec f,\rho)\rho^{n+1}\,d\rho+r^2\int_r^1 \om(\rho)\,M_p(\vec f,\rho)\,d\rho\right)\\
& \leq \begin{cases}
c\,\om(r)\,r\,\| \vec f\|_p  & \hbox{for}\ 0<r<1 \\
c\,\delta\,r^{-n}\,\| \vec f\|_p  & \hbox{for}\ r>1.
\end{cases}
\end{aligned}
\]
We see that $\xi_2\in Y$. The proofs that $\xi_3$ and $\xi_4$ are in $Y$ are quite similar to estimating $T_2$ and $T_3$ above, so we will not give the details. But we can conclude not only that $\xi\in Y$, but
\begin{equation}\label{est:xi}
\|\xi\|_Y\leq c\, (\|\vec f\|_p+\| f_0\|_p).
\end{equation}

Now we let $w\in Y$ be the fixed point of $S$, so $w$ satisfies \eqref{generalvariational3}. We use $w$ to find $\widetilde v^w$  and then \eqref{eqnfor-u0'} to find $u_0'$. Integrating \eqref{eqnfor-u0'} to find $u_0$ (up to a constant) and letting $\widetilde v=\widetilde v^w$, we have
\begin{equation}\label{u*}
u^*(x):=u_0(r)+\widetilde x\cdot\widetilde v(r)+ w(x)
\end{equation}
is a solution of \eqref{eq:HalfspaceProblem1}. Now we want to show that $u^*$ has the desired regularity properties. Since $w=(I+T)^{-1}\xi\in Y$, we know that $M_p(\nabla w,r)\leq c\,\om(r)$ as $r\to 0$. Moreover, $Pw=0$ implies that $\int_{|x|<r} w\,dx=0$ for every $r>0$. Using this and $p>n$, Morrey's inequality (cf.\ \cite{GT}) implies
\[
\sup_{|x|<r}|w(x)| \leq c_n\,r\, \left(\meanint_{|y|<r}|\nabla w|^p\,dy\right)^{1/p}.
\]
(Recall that $|x|<r$ still refers to points $x\in\RR^n_+$.)
But for fixed $r\in (0,1)$ we can introduce $r_j=2^{-j}\,r$ and compute
\[
\meanint_{|y|<r} |\nabla w|^p\,dy=\frac{n}{r^n|S^{n-1}_+|} \sum_{j=0}^\infty \int_{r_{j+1}<|y|<r_j} |\nabla w|^p\,dy
\leq c\sup_{0<\rho<r} M^p_p(\nabla w,\rho).
\]
We conclude that
\begin{equation}
\sup_{|x|<r} |w(x)|\leq c\,r\,\om(r)\quad \hbox{as $r\to 0$,}
\end{equation}
which implies that $w$ is differentiable at $x=0$ with $\partial_jw(0)=0$ for $j=1,\dots,n$.
Moreover, our assumption that \eqref{eq:DynSystm}  is uniformly stable as $t\to\infty$ implies by Proposition \ref{pr:stabilty} in Appendix E that $(\phi,\psi)$ remains bounded as $t\to\infty$, and in fact $|\psi(t)|\leq c\,\e(t)$ as $t\to\infty$. We now want to use \eqref{ourODEsystem:f} to show that $\widetilde v$ is bounded as $t\to\infty$. 
From the second component in 
\eqref{ourODEsystem:f} we have
\[
|\widetilde v_t(t)|\leq  c_1\,\e(t) +c_2\,\meanint|\nabla w|ds.
\]
Let us integrate this from $T$ to $T+\ln 2$:
\[
\int_T^{T+\ln 2} |\widetilde v_t(t)|\,dt \leq  c_1\, \e(T) + c_2\,\int_T^{T+\ln 2} \int_{S^{n-1}_+} |\nabla w|\,ds\,dt.
\]
But letting $R=e^{-T-\ln 2}$ we find
\[
\int_T^{T+\ln 2} \int_{S^{n-1}_+} |\nabla w|\,ds\,dt = \int_R^{2R} \int_{S^{n-1}_+} |\nabla w| \,ds\,\frac{dr}{r}
\leq \frac{c}{R^n}\int_{A_R^+}|\nabla w|\,dx \leq c\,M_{p}(\nabla w,R).
\]
Since we have assumed that $M_p(\nabla w,r)$ is bounded by $\om(r)$ as $r\to 0$,  we have shown
\begin{subequations}
 \begin{equation}\label{est:int-v_t}
\int_T^{T+\ln 2} |\widetilde v_t(t)|\,dt \leq c\,\e(T) \quad\hbox{as}\ T\to \infty.
\end{equation}
Using this in the first component in \eqref{ourODEsystem:f}, we have
 \begin{equation}\label{est:int-v}
\int_T^{T+\ln 2} |\widetilde v(t)|\,dt \leq C \quad\hbox{as}\ T\to \infty.
\end{equation}
\end{subequations}
But now we may use the elementary inequality
 \begin{equation}
\sup_{a\leq t \leq b}|v(t)| \leq c\int_a^b (|v(t)|+|v_t(t)| )\,dt
\end{equation}
to conclude that $|\widetilde v(T)|$ is bounded as $T\to \infty$.
Of course, $\widetilde v(t)$ is actually $\widetilde v(e^{-t})=\widetilde v(r)$, so we see that $|\widetilde v(r)|$ is bounded as $r\to 0$, 
and hence $\widetilde x\cdot\widetilde v$ is Lipschitz.
Finally, using \eqref{eqnfor-u0'}, we can estimate
\[
\begin{aligned}
|u_0(r)-u_0(0)|  \leq \int_0^r|u_0'(\rho)|\,d\rho 
 \leq c\int_0^r \left(|\vartheta(\rho)|+\om(\rho)\rho|\widetilde v\,'(\rho)|+\om(\rho)|\widetilde v(\rho)|
+\om(\rho)\meanint |\nabla w|\,ds\right)\,d\rho.
\end{aligned}
\]
From \eqref{def:vartheta} we have $\vartheta(r)=\overline{f}_1(r)+r^{1-n}\int_0^r\overline{f_0}(\rho)\rho^{n-1}\,d\rho$ where
$\overline{f}_1,\overline{f_0}$ are given in \eqref{def:f1,f0}. Since $f_0\in L^p(\RR^n)$ and vanishes for $r>1$, we can estimate
$
\int_0^\rho |\overline{f_0}(\tau)|\tau^{n-1}\,d\tau 
=c\int_{B_\rho}|f_0(x)|\,dx 
\leq c\,\rho^n\|f_0\|_{L^p}.
$
Hence
\[
\int_0^r\rho^{1-n}\int_0^\rho |\overline{f_0}(\tau)|\tau^{n-1}\,d\tau
\leq c\,\|f_0\|_{L^p}\int_0^r\rho\,d\rho =c\,r^2\, \|f_0\|_{L^p}.
\]
Since $\overline{f}_1$  vanishes for $r>1$ and also for $0<r<1/4$, we can even more easily verify that
$
\int_0^r |\overline{f_1}(\rho)|\,d\rho
\leq c\,r^2\,\|\vec f\|_{L^p}.
$
Since $\widetilde v(r)$ is bounded $r\to 0$, we see that 
$\int_0^r \om(\rho)|\widetilde v(\rho)| \,d\rho\leq c\,r\,\om(r)$.
To estimate the last term, we can proceed similarly to \eqref{int-nabla_w-est} (but letting $r_j=2^{-j}\,r$) to conclude
\[
\begin{aligned}
\int_0^r \om(\rho)\,\meanint_{S^{n-1}_+} |\nabla w|\,ds\,d\rho\leq\sum_{j=0}^\infty \om(r_j)\,r_j\,M_p(\nabla w,r_j) \\
\leq c\,\sum_{j=0}^\infty \om^2(r_j)\,r_j \leq c\,\int_0^r\om^2(\rho)\,d\rho = o(r)\quad\hbox{as $r\to 0$},
\end{aligned}
\]
since $\int_0^r\rho^{-1}\om^2(\rho)\,d\rho\to 0$ as $r\to 0$.
We have one more term to estimate (using $r_j=2^{-j}\,r$):
\[
\begin{aligned}
\int_0^r  \om(\rho) \rho |\widetilde v'(\rho)|\,d\rho
=\sum_{j=0}^\infty \int_{r_{j+1}}^{r_j} \om(\rho)\,\rho\,\widetilde v\,'(\rho)\,d\rho
\leq \sum_{j=0}^\infty \om(r_j) \int_{t_j}^{t_j+\ln 2}\widetilde v_t e^{-t}\,dt \\
\leq \sum_{j=0}^\infty r_j\,\om^2(r_j) \leq c\,\int_0^r \om^2(\rho) \,d\rho= o(r)\quad\hbox{as $r\to 0$},
\end{aligned}
\]
where we have used \eqref{est:int-v_t}).
We conclude that
\begin{equation}
|u_0(r)-u_0(0)|
=o(r) \quad\hbox{as}\ r\to 0,
\end{equation}
which shows that $u_0$ is differentiable at $r=0$ with $u_0'(0)=0$. Since $u_0$ and $w$ in \eqref{u*} are differentiable and 
$\widetilde x\cdot\widetilde v$ is Lipschitz at $x=0$, we conclude that $u^*$ is Lipschitz  at $x=0$.

Next we need to confirm that $u=u^*$ in order to conclude that $u$ is Lipschitz at $x=0$. But $u$ and $u^*$ both satisfy \eqref{eq:HalfspaceProblem1} and the estimate $M_{1,p}(u,r)\leq c\,r^{-n}$ as $r\to \infty$. Then, by Corollary \ref{co:uniqueness} in Appendix F, we see indeed that $u=u^*$.

Finally, let us also assume that all solutions of \eqref{eq:DynSystm}  are asymptotically constant. Then, by 
Proposition \ref{pr:stabilty} in Appendix E, we know that $\phi(t)\to\phi_\infty$ as $t\to\infty$. Using \eqref{ourODEsystem:f},
we can apply the above arguments to $\widetilde v-n\phi_\infty$ to conclude
\[
\sup_{0<\rho<r}|\widetilde v(\rho)-n\phi_\infty|\leq c\,\om(r).
\]
 This shows that $\widetilde x\cdot\widetilde v(r)$ is differentiable at $x=0$. Putting this together with the differentiabilty of $u_0$ and $w$ at $x=0$, we have completed the proof of Theorem \ref{th:1}. \hfill$\Box$
 
 \medskip
We can use the results of \cite{MM3} on  the largest eigenvalue $\mu(r)$ of the 
 of the symmetric matrix $S(r)=-\frac{1}{2}(R(r)+R^t(r))$ to obtain the following corollaries of Theorem \ref{th:1}; in both we assume  the $a_{ij}$ satisfy \eqref{aij-deltaij}, where $\om$ satisfies \eqref{eq:Sq-Dini} and 
 \eqref{om-vanishing}.
 To begin with, in \cite{MM3} it is shown that \eqref{mu-cond1} implies that \eqref{eq:DynSystm} is uniformly stable;
 hence we obtain the following:

\begin{corollary} \label{co:1} 
Suppose that $\mu(r)$ satisfies \eqref{mu-cond1}. Then every  solution $u\in H_\loc^{1,2}(\RR^n_+)$ of \eqref{eq:HalfspaceProblem} is Lipschitz continuous at $x=0$.
 \end{corollary}

\noindent
Moreover,  in \cite{MM3} it is shown that \eqref{mu-cond2} implies that the null solution of \eqref{eq:DynSystm} is asymptotically stable, which in turn shows that $\widetilde v$ in \eqref{spectral-decomp} tends to zero as $r\to 0$.
Consequently, we obtain the following:

\begin{corollary} \label{co:2} 
Suppose that  $\mu(r)$  satisfies \eqref{mu-cond2}. Then every  solution $u\in H_\loc^{1,2}(\RR^n_+)$ of \eqref{eq:HalfspaceProblem} is differentiable at $x=0$ and all derivatives are zero: $\partial_j u(0)=0$ for $j=1,\dots,n$.
 \end{corollary}

\section{Curved Boundaries}

In this section we consider the regularity of a weak solution of \eqref{eq:NeumannProblem} near a point on $\partial U$. Since we are interested in the local behavior of solutions, we may assume $U$ is bounded, the point on $\partial U$ is the origin in $\RR^n$, and the boundary $\partial U$ is given near the origin by $x_n=h(\widetilde x)$ where $h(\widetilde 0)=0$. Recall our assumption
\eqref{h'-sqDini}, which implies that $h$ is differentiable at $\widetilde x=0$ and $\nabla h(\widetilde 0)=0$. 

Let us introduce new independent variables
\[
y_j=x_j \ \hbox{for $j=1,\dots,n-1$} \quad \hbox{and}\quad y_n=x_n-h(x_1,\dots,x_{n-1}).
\]
Notice that $\partial y_j/\partial x_k =\delta_{jk}$ for $j\not= n$ and $\partial y_n/\partial x_k =-\partial h/\partial x_k$ for $k=1,\dots,n-1$ and $\partial y_n/\partial x_n=1$. Consequently, the Jacobian determinant for this change of variables is  1 and by the chain rule
\[
\frac{\partial u}{\partial x_k}=\frac{\partial u}{\partial y_k}-\frac{\partial u}{\partial y_n}\frac{\partial h}{\partial x_k}\ \hbox{for}\ k=1,\dots,n-1 \quad\hbox{and}\quad \frac{\partial u}{\partial x_n}=\frac{\partial u}{\partial y_n}.
\]
We want to express \eqref{eq:weakform} in terms of the $y$-coordinates. Let $i'$ and $j'$ be indices that range from $1$ to $n-1$. 
Then
\[
\begin{aligned}
a_{ij}  \,\frac{\partial u}{\partial x_j}\,\frac{\partial \eta}{\partial x_i}  =&\, a_{i'j'}\left(\frac{\partial u}{\partial y_{j'}}-\frac{\partial u}{\partial y_n}\frac{\partial h}{\partial x_{j'}}\right)
\left(\frac{\partial \eta}{\partial y_{i'}}-\frac{\partial \eta}{\partial y_n}\frac{\partial h}{\partial x_{i'}}\right) \cr
& +
a_{i' n}\frac{\partial u}{\partial y_n}\left(\frac{\partial \eta}{\partial y_{i'}}-\frac{\partial \eta}{\partial y_n}\frac{\partial h}{\partial x_{i'}}\right) 
+ a_{nj'}\left(\frac{\partial u}{\partial y_{j'}}-\frac{\partial u}{\partial y_n}\frac{\partial h}{\partial x_{j'}}\right)\frac{\partial\eta}{\partial y_n}
+ a_{nn}\frac{\partial u}{\partial y_n}\frac{\partial \eta}{\partial y_n}\cr
= &\,a_{i'j'}\frac{\partial u}{\partial y_{j'}}\frac{\partial \eta}{\partial y_{i'}}
+ \left(a_{i'n}-a_{i'j'}\frac{\partial h}{\partial x_{j'}}\right)\frac{\partial u}{\partial y_n}\frac{\partial\eta}{\partial y_{i'}}
+ \left(a_{nj'}-a_{i'j'}\frac{\partial h}{\partial x_{i'}}\right)\frac{\partial u}{\partial y_{j'}}\frac{\partial\eta}{\partial y_{n}}\cr
&+\left( a_{nn}-a_{i'n}\frac{\partial h}{\partial x_{i'}}-a_{nj'}\frac{\partial h}{\partial x_{j'}} +a_{i'j'}\frac{\partial h}{\partial x_{j'}}\frac{\partial h}{\partial x_{i'}}\right)
\frac{\partial u}{\partial y_n}\frac{\partial \eta}{\partial y_n}.
\end{aligned}
\]
Now, if we let $U_0=U\cap B_\e(0)$ for $\e>0$ sufficiently small, then $x\in U_0$ satisfies $x_n>h(x_1,\dots,x_{n-1})$, so if we let $V_0$ denote the corresponding domain in the $y$-variables, then $V_0\subset\RR^n_+$ and
\[
\int_{U_0} a_{ij}\frac{\partial u}{\partial x_i}\frac{\partial \eta}{\partial x_j}\,dx=
\int_{V_0} \widetilde a_{ij}\frac{\partial u}{\partial y_i}\frac{\partial \eta}{\partial y_j}\,dy,
\]
where
\[
\widetilde a_{ij}=\begin{cases}
a_{ij} & {\rm if}\ 1\leq i,j \leq n-1, \\
a_{in}-a_{ij'}\frac{\partial h}{\partial x_{j'}} & {\rm if}\ 1\leq i\leq n-1, \ j=n, \\
a_{nj}-a_{i'j}\frac{\partial h}{\partial x_{i'}} & {\rm if}\ 1\leq j\leq n-1, \ i=n, \\
 a_{nn}-a_{i'n}\frac{\partial h}{\partial x_{i'}}-a_{nj'}\frac{\partial h}{\partial x_{j'}} +a_{i'j'}\frac{\partial h}{\partial x_{j'}}\frac{\partial h}{\partial x_{i'}} & {\rm if}\ i=j=n.
\end{cases}
\]
This enables us to consider the original problem as one for  the coefficients $\widetilde a_{ij}$ in the half-space $\{(y_1,\dots,y_n):y_n>0\}$. In order to apply our results from the previous section, we need $\widetilde a_{ij}$ to be square-Dini continuous and satisfy \eqref{aij-deltaij}; but these conditions follows from our assumption \eqref{h'-sqDini}.

 Now we can write down the 1st-order  dynamical system \eqref{ourODEsystem:joint} associated with the $\widetilde a_{ij}$ in $\RR^{n-1}_+$ whose stability properties determine the differentiability of a weak solution. In particular, the formula \eqref{def:R-halfspace} for the
 $(n-1)\times(n-1)$ matrix $R$ yields
 \begin{equation}\label{R-curved}
 \left[R(r)\right]_{\ell k} =\meanint_{S_+^{n-1}} \left( a_{\ell k} -n\sum_{j=1}^n a_{\ell j}\theta_j\theta_k
 +n\sum_{j=1}^{n-1} a_{\ell j} \frac{\partial h}{\partial x_j}\theta_n\theta_k\right)\,ds_\theta.
 \end{equation}
 In \eqref{R-curved} we need to emphasize that the integrand is considered as a function of $y\in\RR^n_+$, even though the coefficients $a_{ij}$ and $h$ were originally defined in the $x$ variables. Also, note that if $h\equiv 0$, then we are in the half-space case, and the formula for $R(r)$ in \eqref{R-curved} agrees with \eqref{def:R-halfspace}.
 
 \begin{theorem} 
Suppose that $U$ is a bounded domain with Lipschitz boundary $\partial U$ containing the point $0$, near which the boundary can be represented as $x_n=h(x_1,\dots,x_{n-1})$. Suppose  the
$a_{ij}$ satisfy \eqref{aij-deltaij} and $h$ satisfies   \eqref{h'-sqDini}, where $\om$ satisfies \eqref{eq:Sq-Dini} and  \eqref{om-vanishing}. If  the dynamical system \eqref{eq:DynSystm} with matrix $R$ given by \eqref{R-curved} is uniformly stable as $t\to\infty$, then every solution $u\in H^{1,2}(U\cap B)$ of \eqref{eq:weakform} is Lipschitz continuous at $x=0$. If, in addition, every solution of the dynamical system 
\eqref{eq:DynSystm}, \eqref{R-curved} is asymptotically constant as $t\to\infty$, then $u$ is differentiable at $x=0$.
\end{theorem}

As in Section 2, the conditions \eqref{mu-cond1} and \eqref{mu-cond2} can be used to obtain corollaries of this theorem, but now $\mu$ is the largest eigenvalue of the matrix $S(r)=-\frac{1}{2}(R(r)+R^t(r))$, where 
$R$ is given by \eqref{R-curved}. In the following results we assume the conditions on $a_{ij}$, $h$, and $\om$ stated in the theorem; in the second one we note that $\partial u(0)/\partial x_j=\partial u(0)/\partial y_j$
since $\nabla h(\widetilde 0)=0$.

\begin{corollary} \label{co:3} 
Suppose that $\mu(r)$ satisfies \eqref{mu-cond1}. Then every  solution $u\in H_\loc^{1,2}(U\cap B)$ of \eqref{eq:weakform} is Lipschitz continuous at $x=0$.
 \end{corollary}

\begin{corollary} \label{co:4} 
Suppose that  $\mu(r)$  satisfies \eqref{mu-cond2}. Then every  solution $u\in H_\loc^{1,2}(U\cap B)$ of \eqref{eq:weakform} is differentiable at $x=0$ and all derivatives are zero: $\partial_j u(0)=0$ for $j=1,\dots,n$.
 \end{corollary}

  \section{Examples: $n=2$}
  
Let us first consider variable coefficients $a_{ij}$ in $\RR^2_+$. For $n=2$ we have $\theta_1=\cos\phi$ and $\theta_2=\sin\phi$ for $0<\phi<\pi$, so
\eqref{def:R-halfspace} yields a scalar function
\begin{subequations}
\begin{equation}\label{def:R-dim2}
R(r) = \frac{1}{\pi}\int_0^\pi \left(a_{11}(r\theta)-2a_{11}(r\theta)\cos^2\phi
-2a_{12}(r\theta)\cos\phi\sin\phi\right)d\phi.
\end{equation}
In this case, the dynamical system (\ref{eq:DynSystm}) is just a single equation, and we easily find the general solution:
$\phi(t)=C\, \exp[-\int_T^t R(e^{-\tau})\,d\tau]$. Consequently (cf.\ Remark \ref{rk:2} in Appendix D), we know that (\ref{eq:DynSystm}) is uniformly stable if and only if 
$\int_s^t R(e^{-\tau})\,d\tau$ is uniformly bounded below for $T<s<t<\infty$. Expressing this  in terms of $r$ rather than $t$, we see that uniform stability
\begin{equation}\label{unifstable-dim2}
\int_{r_1}^{r_2}\frac{R(\rho)}{\rho}\,d\rho > -K \ \hbox{ for all $0<r_1<r_2<\e$}
\end{equation}
implies that every weak solution of \eqref{eq:HalfspaceProblem} is Lipschitz at the origin.
Similarly, solutions of (\ref{eq:DynSystm}) are asymptotically constant when $\int_T^\infty R(e^{-\tau})\,d\tau$ either converges to a finite number or diverges to $\infty$. In terms of $R(r)$, we find that \eqref{unifstable-dim2} together with
\begin{equation}\label{asymconst-dim2}
\int_0^\e \frac{R(\rho)}{\rho}\,d\rho  \ \hbox{converges to an extended real number $>-\infty$}
\end{equation}
\end{subequations}
imply that every weak solution of \eqref{eq:HalfspaceProblem} is differentiable at  the origin.

To make all this more precise, let us turn to a class of operators considered in  \cite{GS} and \cite{MM3}:
\begin{subequations}
\begin{equation}\label{GS}
a_{ij}=\delta_{ij}+g(r)\theta_i\theta_j\,,
\end{equation}
where $|g(r)|\leq c\,\om(r)$. In this case we can calculate $R(r)=-\frac{1}{2}\,g(r)$ so that uniform stability
\begin{equation}\label{unifstable-GS}
\int_{r_1}^{r_2}\frac{g(\rho)}{\rho}\,d\rho < K \ \hbox{for all $0<r_1<r_2<\e$}
\end{equation}
implies that every weak solution of \eqref{eq:HalfspaceProblem} is Lipschitz continuous at the origin; and if in addition
\begin{equation}\label{asymconst-GS}
\int_0^\e \frac{g(\rho)}{\rho}\,d\rho  \ \hbox{converges to an extended real number $<\infty$,}
\end{equation}
\end{subequations}
then every weak solution of \eqref{eq:HalfspaceProblem} is differentiable at  the origin.

For \eqref{GS} we can  construct explicit solutions of \eqref{eq:HalfspaceProblem} by solving an ODE. For example, if we let
\begin{equation}
u(r,\phi)=U(r)\,\cos\phi,
\end{equation}
then this is a solution provided $U$ satisfies
\begin{equation}
\frac{1}{r}\left[ (1+g(r))\,r\,U'\right]'-\frac{1}{r^2}\,U=0.
\end{equation}
Moreover, we can determine the behavior of $U(r)$ as $r\to 0$ from that of $g(r)$. To do this, it is simpler to again use the variable $t=-\log r$. Letting $\tilde g(t)=g(e^{-t})$, we want $U$ to satisfy
\begin{equation}\label{U-equation}
\frac{d}{dt} \left[ (1+\tilde g(t))\,\frac{dU}{dt}\right]-U=0 \quad \hbox{as}\ t\to\infty.
\end{equation}
We can apply standard results in the asymptotic theory of ODEs. For example, if $\tilde g(t)$ is $C^1$ and satisfies
\begin{equation}\label{g-conditions}
\tilde g(t),\ \frac{d\tilde g}{dt} = o(1) \quad\hbox{as}\ t\to\infty,
\end{equation}
then we can apply Theorem 2.2.1 in \cite{E} to conclude that a solution $U(t)$ of \eqref{U-equation} exists for which both 
$U(t)$ and $(1+\tilde g(t))dU/dt$ are asymptotic to
\begin{equation}
 (1+\tilde g(t))^{-1/4} \exp\left(-\int_1^t \left( \frac{1}{1+\tilde g(s)}+\frac{(d\tilde g/ds)^2}{16(1+\tilde g(s))^2}\right)^{1/2}\,ds\right) 
 \sim e^{-t} \exp\left(\frac{1}{2}\int_1^t \tilde g(s)\,ds\right).
  \end{equation}
  This solution satisfies the finite-energy condition $\int_1^\infty(U^2+(U_t)^2)e^{-nt}\,dt<\infty$, so $u$ is an $H^{1,2}$-solution of \eqref{eq:HalfspaceProblem}. However, if $g(r)$ does not satisfy the Dini condition at $r=0$  then 
 $\int_1^t\tilde g(s)\,ds\to\infty$ as $t\to\infty$ and 
 $u$ is not Lipschitz continuous at  the origin. An example of such a function $g(r)$ is 
  \begin{equation}
  g(r)= |\log r|^{-\alpha} \quad\hbox{where}\ 1/2<\alpha\leq 1;
  \end{equation}
 note that $\tilde g(t)=t^{-\alpha}$  satisfies \eqref{g-conditions} but \eqref{unifstable-GS} is not satisfied.
  In particular, this example shows that a weak solution of \eqref{eq:HalfspaceProblem} when the coefficients $a_{ij}$
  are square-Dini continuous need not be Lipschitz continuous if the associated dynamical system   (\ref{eq:DynSystm})  is not uniformly stable.

Next let us  suppose that the origin lies on the boundary $\partial U$, which locally has the form $x_2=h(x_1)$, where $h(0)=0$ and $|h'(r)|\leq c\,\om(r)$ as $r\to 0$. Then we introduce new independent variables $y_1=x_1$ and $y_2=x_2-h(x_1)$ and consider $a_{ij}$ as functions of $(y_1,y_2)\in\RR^2_+$. We can calculate the scalar function $R(r)$ in \eqref{R-curved}:
  \begin{equation}\label{R-curved,n=2}
  \frac{1}{\pi} \int_0^\pi \left( a_{11}(r\theta)-2\left(a_{11}(r\theta)\cos^2\phi+a_{12}(r\theta)\cos\phi\sin\phi\right)+2a_{11}(r\theta)h'(r\theta_1)\cos\phi\sin\phi\right)\,d\phi.
  \end{equation}
  Again, we find that \eqref{unifstable-dim2} implies that every weak solution $u\in H^{1,2}(U)$ of \eqref{eq:NeumannProblem} is Lipschitz at  the origin, and if \eqref{asymconst-dim2} also holds then $u$ is differentiable there.
  
Now let us consider the special case of \eqref{R-curved,n=2} when the operator is the Laplacian, so that $a_{ij}=\delta_{ij}$. In this case, we have simply
   \begin{equation}\label{R-curved,n=2,Laplace}
  R(r)  =\frac{2}{\pi}\int_0^\pi h'(r\cos\phi)\cos\phi\sin\phi\,d\phi.
  \end{equation}
 One way to make sure that \eqref{unifstable-dim2} and \eqref{asymconst-dim2} hold is to have $R(r)\geq 0$ for $0<r<\e$. This will be the case, for example, if 
 \begin{equation}\label{h-convex}
h'(x) \leq 0 \quad \hbox{for $-\e<x<0$} \quad\hbox{and}\quad h'(x)\geq 0 \quad  \hbox{for $0<x<\e$.}
\end{equation}
Consequently, if the boundary function $h$ satisfies \eqref{h-convex}, we can conclude that every weak solution 
$u\in H^{1,2}(U)$ of \eqref{eq:NeumannProblem} is differentiable at  the origin.
  
  We should compare our results for the Laplacian with those of \cite{W} concerning conformal maps. In \cite{W}, the hypotheses on the boundary are weaker than ours, and asymptotics are obtained, not just conclusions about differentiability. However, under the hypotheses on the boundary that we consider, Theorem XI(A) in \cite{W} shows that the  behavior of a conformal map as $z\to 0$ is dominated by
\begin{equation}\label{W-asym1}
\exp\left[ -\pi\int_{|z|}^a \frac{1}{r\,\Theta(r)}\,dr\right].
\end{equation}
Here $\Theta(r)$ measures the angle between the two arcs $\Gamma_-$ and $\Gamma_+$ corresponding to $x_2=h(x_1)$ for $x<0$ and for $x>0$ respectively. Consequently, $|\Theta(r)-\pi|\leq \om(r)$ as $r\to 0$, and we can write
\[
\frac{\pi}{r\Theta(r)}=\frac{1}{r}\left[1-\left(1-\frac{\Theta(\rho)}{\pi}\right) \right]^{-1}
=\frac{1}{r}\left[ 1+ \left(1-\frac{\Theta(r)}{\pi}\right) + O(\om^2(r)) \right].
\]
Thus, as $|z|\to 0$, \eqref{W-asym1} is asymptotic to
\begin{equation}\label{W-asym2}
 C\,|z|\,\exp\left[\frac{1}{\pi}\int_{|z|}^a \frac{\Theta(r)-\pi}{r}\,dr \right].
\end{equation}
This means, for example, that the convergence (or divergence to $-\infty$) of the intergal in \eqref{W-asym2} determines whether the conformal map is Lipschitz continuous at $z=0$; this is the analogue to our condition \eqref{unifstable-dim2} for a harmonic function to be Lipschitz continuous at the origin.
    
\appendix
\section{Asymptotic expansion of the Neumann function}
\label{DerivationN(x,y)}

In this appendix we derive the asymptotic expansion of the Neumann function $N(x,y)$. 
 We need to use an expansion of the fundamental solution $\Gamma$ in spherical harmonics. 
Let ${\mathcal H}(k)$ denote the spherical harmonics of degree $k$ and let $N(k)=\hbox{dim}\,{\mathcal H}(k)$. For each $k$, choose a basis $\{\varphi_{k,m}:m=1,\dots,N(k)\}$ for ${\mathcal H}(k)$ that is orthonormal with respect to the spherical mean inner product:
\[
\meanint_{S^{n-1}}\varphi_{k,\ell}\,\varphi_{k,m}\,ds=
\begin{cases} 1 & \ell=k \\ 0 & \ell\not= m \end{cases}.
\]
 For notational convenience, let 
$\hat x=x/|x|$ and $\hat y=y/|y|$. We also assume $n\geq 3$, the case $n=2$ being analogous.
For $|x|<|y|$ we can write
$\Gamma(|x-y|)$ as a convergent series
\begin{subequations}\label{Gamma-expansion:joint}
\begin{equation} \label{Gamma-expansion1}
\Gamma(|x-y|)=\sum_{k=0}^\infty \,\frac{|x|^k}{|y|^{n-2+k}}\sum_{m=1}^{N(k)}a_{k,m}\,\varphi_{k,m}\left(\hat{x}\right)\,\varphi_{k,m}\left(\hat{y}\right),
\end{equation}
where $a_{k,m}$ are  certain coefficients.\footnote{The expansion \eqref{Gamma-expansion1} was used in \cite{MM3}, but the coefficients were unfortunately left out of the formula there.}  With $x=0$ we know that $\Gamma(|y|)=a_0\,|y|^{2-n}$ with $a_0=(2-n)^{-1}\om_n^{-1}$ where $\om_n=|S^{n-1}|$. We can also use a Taylor series for $f_y(x)=|x-y|^{2-n}$, i.e.\ 
\[
|x-y|^{2-n}=|y|^{2-n}+(n-2)|y|^{-n}\sum_j x_j y_j +\cdots
\]
 to compute the other coefficients. For example, we can write
 \begin{equation}\label{Gamma-expansion2}
 \Gamma(|x-y|)=a_0\left(\frac{1}{|y|^{n-2}}+(n-2)\frac{|x|}{|y|^{n-1}} \sum_{m=1}^n \hat x_m\,\hat y_m+\cdots\right).
 \end{equation}
 \end{subequations}
But to compute  our Neumann function $N(x,y)$, we want the basis $\{\varphi_{km}\}$ for $k>1$ to also possess certain symmetries with respect to the half-space.

Recall that the spherical harmonics of degree $k$ are generated by the restriction to the unit sphere of the harmonic polynomials of degree $k$:
\begin{equation}
h(x)=\sum_{|\alpha|=k} c_\alpha x^\alpha\quad\hbox{is harmonic},
\end{equation}
where $\alpha=(\alpha_1,\dots,\alpha_n)$ and $x^\alpha=x_1^{\alpha_1}\cdots x_n^{\alpha_n}$.
For our half-space geometry, we want to distinguish those harmonic functions for which $\alpha_n$ is even or odd. 
Let ${\calH}_{e}(k)$ be the spherical harmonics corresponding to even $\alpha_n$ 
 and let $N_e(k)$ denote its dimension; choose an orthonormal basis $\{\varphi^{e}_{k,m}: m=1,\dots, N_e(k)\}$ for $\calH_e(k)$. Similarly, let $\calH_{o}(k)$ be the spherical harmonics corresponding to odd $\alpha_n$ and choose an orthonormal basis
$\{\varphi^{o}_{k,m}: m=1,\dots, N_o(k)\}$ for  $\calH_{o}(k)$. Then 
$\{\varphi^{e}_{k,m}: m=1,\dots, N_e(k)\}\cup \{\varphi^{o}_{k,m}: m=1,\dots, N_o(k)\}$ is an orthonormal basis for ${\calH}(k)$
which we may use to rewrite \eqref{Gamma-expansion1} as:
\begin{equation} \label{eq:Gamma-expansion2}
\Gamma(|x-y|)=\sum_{k=0}^\infty \,\frac{|x|^k}{|y|^{n-2+k}}
\left(\sum_{m=1}^{N_e(k)}a_{k,m}\,\varphi^e_{k,m}\left(\hat{x}\right)\,\varphi^e_{k,m}\left(\hat{y}\right) +
\sum_{m=1}^{N_o(k)}b_{k,m}\,\varphi^o_{k,m}\left(\hat{x}\right)\,\varphi^o_{k,m}\left(\hat{y}\right)
\right),
\end{equation}
But 
\[
y^*=(\tilde y,-y_n) \quad\Rightarrow\quad
\varphi^e_{k,m}(\hat y^*)=\varphi^e_{k,m}(\hat y) \quad \hbox{and} \quad \varphi^o_{k,m}(\hat y^*)=-\varphi^o_{k,m}(\hat y)
\]
so
\[
\Gamma(|x-y^*|)=\sum_{k=0}^\infty \,\frac{|x|^k}{|y|^{n-2+k}}
\left(\sum_{m=1}^{N_e(k)}a_{k,m}\,\varphi^e_{k,m}\left(\hat{x}\right)\,\varphi^e_{k,m}\left(\hat{y}\right) -
\sum_{m=1}^{N_o(k)}b_{k,m}\,\varphi^o_{k,m}\left(\hat{x}\right)\,\varphi^o_{k,m}\left(\hat{y}\right)\right).
\]
When we add $\Gamma(|x-y|)$ and $\Gamma(|x-y^*|)$ the terms involving $\varphi^o_{k,m}$ cancel, so we obtain
\[
N(x,y)=2\sum_{k=0}^\infty \,\frac{|x|^k}{|y|^{n-2+k}}
\sum_{m=1}^{N_e(k)}a_{k,m}\,\varphi^e_{k,m}\left(\hat{x}\right)\,\varphi^e_{k,m}\left(\hat{y}\right) 
\quad\hbox{for $|x|<|y|$}.
\]

Restricting the $\varphi^e_{k,m}$ to $S^{n-1}_+$ yields spherical harmonics  with zero normal derivative along the boundary $\partial S^{n-1}_+$, but we also want them to be orthonormal with respect to the spherical mean inner product on $S^{n-1}_+$. We easily calculate
\[
\meanint_{S_+^{n-1}}\varphi^e_{k,m}\varphi^e_{k',m'}\,ds=
\frac{1}{2}\,\meanint_{S^{n-1}}\varphi^e_{k,m}\varphi^e_{k',m'}\,ds=
\begin{cases}\frac{1}{2} & k=k' \ \hbox{and}\  m=m', \\ 0 & k\not=k' \ \hbox{or}\ m\not=m'.
\end{cases}
\]
Consequently, we will have an orthonormal basis $\{\tilde\varphi_{k,m}:m=1,\dots,\tilde N(k)\}$ of spherical harmonics with zero normal derivative along the boundary $\partial S^{n-1}_+$ if we define: 
\begin{equation}\label{def:tilde-varphi}
\tilde\varphi_{k,m}=\sqrt{2}\,\varphi^e_{k,m}|_{S^{n-1}_+} \quad\hbox{and}\quad \tilde N(k)=N_e(k).
\end{equation}
For $k=1$, we want  $\tilde \varphi_{1,m}=\bar c\, \theta_m$ for some constant $\bar c$ and all $m=1,\dots,n-1$.
Using \eqref{def:c_n} and the fact that the $\varphi_{1,m}$ are orthonormal, we see that 
\begin{equation}\label{eq:phi_{1,m}}
\tilde\varphi_{1,m}=\frac{1}{\sqrt{c_n}}\,\theta_m, \quad\hbox{for $m=1,\dots,n-1$.}
\end{equation}
 We  therefore obtain \eqref{eq:N-x<y}.
By interchanging the roles of $x$ and $y$ we  get the expansions of $\Gamma(|x-y|)$ and $\Gamma(|x-y^*|)$ for $|x|>|y|$, and add them together to obtain \eqref{eq:N-y<x}.

\section{Orthogonality Properties}
\label{OrthogProps}

In this appendix we discuss orthogonality properties necessary to show \eqref{grad(u)-props}. In fact, we first prove the following:
\begin{lemma} 
If $f\in H^{1,1}_{\loc}(\overline{\RR^n}\backslash\{0\})$ and $r>0$, then for $i=1,\dots,n$ we have
\begin{equation}\label{L1}
\meanint_{S^{n-1}_+} f(r\theta)\,ds =0 \ \Rightarrow \ \meanint_{S^{n-1}_+} \theta_i \partial_i f(r\theta) ds=0,
\end{equation}
and for $j=1,\dots,n-1$ we have
\begin{equation}\label{L2}
\meanint_{S^{n-1}_+} \theta_j f(r\theta)\,ds =0 \ \Rightarrow \ \meanint_{S^{n-1}_+}  \partial_j f(r\theta) ds = 0
= \meanint_{S^{n-1}_+}  \theta_j\theta_i\partial_i f(r\theta) ds.
\end{equation}
\end{lemma}  

\noindent{\bf Proof.}
To prove \eqref{L1} we consider $\phi\in C_{comp}^\infty(0,\infty)$ and write
\[
\biggl\langle \meanint \theta_i\del_i f ds,\phi\biggr\rangle
=\int_0^\infty\meanint_{S^{n-1}_+}\,\theta_i\,\del_if(r\theta)ds\,\phi(r)\,dr
=\frac{1}{|S_+^{n-1}|}\int_{\RR_+^n} x_i\del_i f(x)\phi(|x|)|x|^{-n}\,dx.
\]
Taking the divergence of $x_i f(x) \phi(|x|) |x|^{-n}$, we obtain
\[
\partial_i(x_i f(x) \phi(|x|) |x|^{-n})=\theta_i\partial_i f(x) \phi(|x|) |x|^{-n+1} +  f(x) \phi'(|x|) |x|^{-n}.
\]
By the divergence theorem,
\[
\int_{\RR_+^n} \partial_i(x_i f(x) \phi(|x|) |x|^{-n}) dx=-\int_{\RR^{n-1}} (x_n f(x)\phi(|x|) |x|^{-n})  |_{x_n=0} \,d\widetilde x=0,
\]
so 
\[
\int_{\RR^n_+} x_i\partial_i f(x) \phi(|x|) |x|^{-n} dx=-\int_{\RR^n_+} f(x) \phi'(|x|) |x|^{-n+1} dx
=-\int_0^\infty \int_{S^{n-1}_+} f(r\theta)\,ds \,\phi'(|x|) \,dr.
\]
Using the hypothesis in \eqref{L1}, this last integral vanishes, which confirms the conclusion in \eqref{L1}.

To prove \eqref{L2}, we again consider $\phi\in C_{comp}^\infty(0,\infty)$ and write
\[
\biggl\langle \meanint \,\del_j f ds,\phi\biggr\rangle
=\int_0^\infty \meanint_{S_+^{n-1}}\del_j f(r\theta)\,ds\,\phi(r)\,dr
=\frac{1}{|S_+^{n-1}|} \int_{\RR_+^n}\del_jf(x)|x|^{1-n}\phi(|x|)\,dx.
\]
But  we can integrate by parts in this last integral to obtain
\[
 \int_{\RR_+^n} f(x)[r^{1-n}\phi(r)]'|_{r=|x|}\theta_j\,dx
=\int_0^\infty\int_{S_+^{n-1}}f(r\theta)\theta_j\,ds\,[r^{1-n}\phi(r)]' r^{n-1}dr.
\]
This gives the first conclusion in (\ref{L2}). To obtain the second conclusion, we write
\[
\biggl\langle \meanint \,\theta_j\theta_i\del_i f ds,\phi\biggr\rangle=
\frac{1}{|S_+^{n-1}|}\int_{\RR_+^n} x_j\,x_i\del_i f(x)\phi(|x|)|x|^{-n-1}\,dx.
\]
Take the divergence (for fixed $j$):
\[
\begin{aligned}
\partial_i(x_j\,x_i\,f(x)\phi(|x|)|x|^{-n-1})=
x_j f(x)\phi(|x|)|x|^{-n-1}+ nx_j f(x)\phi(|x|)|^{-n-1} \\
+x_jx_i\partial_i f(x) \phi(|x|) |x|^{-n-1}+x_jx_if(x)\left(\phi(r)r^{-n-1}\right)'\theta_i.
\end{aligned}
\]
So applying the divergence theorem yields
\[
\begin{aligned}
\int_{\RR_+^n} x_j\,x_i\del_i f(x)\phi(|x|)|x|^{-n-1}\,dx=-\int_{\RR^{n-1}}\left(x_jx_nf(x)\phi(|x|)|x|^{-n-1}\right)|_{x_n=0} \,d\widetilde x \\
-\int_{\RR^n_+} \left((n+1)x_jf(x)\phi(|x|)|x|^{-n-1} +x_jx_i f(x) \left(\phi(r)r^{-n-1}\right)'\theta_i \right)dx.
\end{aligned}
\]
The boundary integral clearly vanishes and the domain integral simplifies considerably to yield
\[
\int_{\RR_+^n} x_j\,x_i\del_i f(x)\phi(|x|)|x|^{-n-1}\,dx=-\int_0^\infty\int_{S^{n-1}_+} \theta_j f(r\theta)\,ds\,\phi'(r)\,dr
\]
Using the hypothesis in \eqref{L2}, this last integral vanishes, which confirms the second conclusion in \eqref{L2}.
\hfill$\Box$

Now we are able to address \eqref{grad(u)-props}.

\begin{corollary} 
If $u\in H^{1,2}(B_+(1))$ and we introduce the spectral decomposition \eqref{spectral-decomp}, then there is a constant $c>0$ such that
\[
\int_{B_+(1)} |\nabla u|^2\,dx \geq c \int_0^1 \left[ (u'_0)^2+|\widetilde v|^2 +r^2|\widetilde v\,'|^2\right] r^{n-1}\,dr
+\int_{B_+(1)} |\nabla w|^2\,dx.
\]
\end{corollary}

\noindent{\bf Proof.} We compute 
\[
\nabla_i u(x)= \begin{cases}
u_0'(r)\theta_i + \widetilde v\,'(r) \cdot \widetilde x \ \theta_i + v_i(r) + \nabla_i w, & 1\leq i \leq n-1 \\
u_0'(r)\theta_n + \widetilde v\,'(r) \cdot \widetilde x \ \theta_n  + \nabla_n w, & i=n,
\end{cases}
\]
and
\[
\begin{aligned}
|\nabla u|^2 =& (u_0')^2 + 2 u_0'(\widetilde v\,'\cdot\widetilde x) +2 u_0'(\widetilde\theta\cdot\widetilde v) +
2u_0'(\theta\cdot\nabla w) + (\widetilde v\,' \cdot \widetilde x)^2 +2(\widetilde v\,'\cdot\widetilde x)(\widetilde\theta\cdot\widetilde v)\\
&+2(\widetilde v\,'\cdot\widetilde x)(\theta\cdot\nabla w)+|\widetilde v|^2 + 2\widetilde v\cdot\nabla w+|\nabla w|^2.
\end{aligned}
\]
(In the formula for $|\nabla u|^2$, note that $\widetilde\theta=(\theta_1,\dots,\theta_{n-1})$ and dot products involving  $\widetilde\theta$ or $\widetilde v$ are summed only over $1,\dots,n-1$.) The integral over $S_+^{n-1}$ of some of these terms vanish due to $\int_{S_+^{n-1}}\theta_i\,ds=0$ for $i=1,\dots,n-1$:
\[
\int_{S_+^{n-1}} u'_0\, (\widetilde v \,'\cdot \widetilde x) \,ds=0=\int_{S_+^{n-1}} u'_0\,( \widetilde\theta\cdot\widetilde v) \,ds.
\]
Other terms vanish utilizing \eqref{L1} and \eqref{L2}:
\[
\int_{S^{n-1}_+} u_0'(\theta\cdot\nabla w) ds=0 = \int_{S^{n-1}_+} \widetilde v\cdot\nabla w\, ds
= \int_{S^{n-1}_+}(\widetilde v\,'\cdot\widetilde x)(\theta\cdot\nabla w)\,ds.
\]
Still other terms simplify using \eqref{def:c_n}:
\[
\int_{S^{n-1}_+} (\widetilde v\,'\cdot \widetilde x)^2 ds=
r^2\sum_{i,j=1}^{n-1} v_i' \,v_j' \int_{S^{n-1}_+}\theta_i\theta_j ds=\frac{r^2|S_+^{n-1}|}{n}\sum_{i=1}^{n-1} (v'_i)^2
\]
\[
\int_{S^{n-1}_+} (\widetilde v'\cdot\widetilde x)(\widetilde\theta\cdot\widetilde v)ds
=r\sum_{i,j=1}^{n-1} v'_i v_j \int_{S^{n-1}_+}\theta_i\theta_j\,ds = \frac{r|S_+^{n-1}|}{n}\sum_{i=1}^{n-1} v'_i v_i.
\]
 So
\[
\int_{B_+(1)} |\nabla u|^2\,dx = |S^{n-1}_+| \int_0^1 \left( (u'_0)^2 +\frac{r^2}{n}|\widetilde v'|^2 + \frac{2r}{n} \widetilde v'\cdot \widetilde v
+|\widetilde v|^2\right) r^{n-1}\,dr + \int_{B_+(1)} |\nabla w|^2\,dx.
\]
Using
\[
\frac{2r}{n}\widetilde v\cdot \widetilde v\,' = \frac{2r}{n^{2/3}n^{1/3}}\widetilde v\cdot \widetilde v\,' \geq -\frac{r^2 |\widetilde v\,'|^2}{n^{4/3}}
-\frac{|\widetilde v|^2}{n^{2/3}}
\]
we find
\[
\begin{aligned}
\int_{|x|<1} |\nabla u|^2\,dx \geq |S^{n-1}| \int_0^1 \left( (u_0')^2 + r^2\left(\frac{1}{n}-\frac{1}{n^{4/3}}\right)|\widetilde v\,'|^2+
\left(1-\frac{1}{n^{2/3}}\right) |\widetilde v|^2 \right) r^{n-1}\,dr
\\
 + \int_{|x|<1}|\nabla w|^2\,dx.
\end{aligned}
\]
Since $n^{4/3}>n$ and $n^{2/3}>1$, this completes the proof. \hfill$\Box$

\section{Sobolev regularity of weak solutions}
\label{SobolevReg}

In this appendix we show that, if  $u\in H_{\loc}^{1,2}(\overline{\RR^n_+})$ satisfies \eqref{eq:HalfspaceProblem}
where the $a_{ij}$ are continuous functions then $u\in H_{\loc}^{1,p}(\overline{\RR^n_+})$ for any $p>2$.
Let us introduce
 the operator ${\mathcal L}$, which is defined on  $v\in H^{1,q}_{\loc}(\overline{\RR^n_+})$ for any $q>1$, and assigns a functional on $H_{\rm comp}^{1,q'}(\overline{\RR^n_+})$ defined by
\begin{equation}\label{def:L}
\langle {\mathcal L}v,\eta\rangle = - \int_{\RR^n_+} a_{ij}\del_jv\,\del_i\eta\,dx
\quad\hbox{for all}\ \eta\in H_{\rm comp}^{1,q'}(\overline{\RR^n_+}).
\end{equation}
In this context, we have assumed that $u\in H_{\loc}^{1,2}(\overline{\RR^n_+})$ is a solution of ${\mathcal L}u=0$,  and we want to conclude that $u\in H_{\loc}^{1,p}(\overline{\RR^n_+})$ for any $p>2$.

The assertion $ u\in H_{\loc}^{1,p}(\overline{\RR^n_+})$ is proved by localizing near a point in $\overline{\RR^n_+}$. Since the issue is on the boundary,
we assume the point is $0$, so it suffices to show that $\phi_0u\in H^{1,p}(B_+)$ for some $\phi_0\in C_0^\infty(B)$ with $\phi_0\equiv 1$ near $0$; here $B=\{x\in \RR^n: |x|<1\}$ and $B_+=\{x\in {\RR^n_+}: |x|<1\}$. By continuity, for any $\e>0$ we can find a $\delta>0$ so that 
$\sup_{|x|\leq \delta} |a_{ij}(x)-\delta_{ij}| \leq \e$. However, for notational convenience we simply assume  a small oscillation condition in $B_+$:
\begin{equation}\label{aij-deltaij<epsilon}
\sup_{|x|\leq 1} |a_{ij}(x)-\delta_{ij}| \leq \e,
 \end{equation}

Let us denote by ${\mathcal L}_0$ the operator \eqref{def:L} with $a_{ij}=\de_{ij}$.
Let $N(x,y)$ be the Neumann function for the Laplacian on $\RR^n_+$, and denote the associated integral operator by ${\mathcal N}$.  Note that for $u\in C^1_{comp}(\overline{\RR^{n}_+})$ we have by Green's identities
 \[
 {\mathcal N}{\mathcal L}_0 u(x)=-\int_{\RR^n_+} \nabla_y N(x,y)\cdot\nabla u(y)\,dy
 =\int_{\RR^n_+} \Lap_y N(x,y) u(y)\,dy =u(x).
 \]
Since any $u\in H^{1,p}_{comp}(\overline{\RR^n_+})$ can be approxmiated by $u_j\in C^1_{comp}(\overline{\RR^{n}_+})$,  we conclude that $ {\mathcal N}{\mathcal L}_0$  is the identity on $H^{1,p}_{comp}(\overline{\RR^n_+})$.

For $\phi_1\in C_{comp}^\infty(B)$ satisfying $\phi_0\,\phi_1=\phi_0$ on $B$, let us write
\[
{\mathcal L}_0(\phi_0 u) + ({\mathcal L}-{\mathcal L}_0)(\phi_0 u)=
 \phi_0 {\mathcal L}u + [{\mathcal L},\phi_0](\phi_1 u)=[{\mathcal L},\phi_0](\phi_1 u),
\]
where we have used \eqref{eq:HalfspaceProblem} to conclude $ \phi_0 {\mathcal L}u=0$.
Now we apply ${\mathcal N}$ to conclude
\begin{equation}\label{I+S}
\phi_0 u + {\mathcal N}({\mathcal L}-{\mathcal L}_0)(\phi_0 u)=
{\mathcal N} [{\mathcal L},\phi_0](\phi_1 u).
\end{equation}
Taking $\e=\e(p)$ sufficiently small in \eqref{aij-deltaij<epsilon}, we can arrange that both
\begin{equation}\label{2_maps_on_B_+}
{\mathcal N}({\mathcal L}-{\mathcal L}_0)\phi_1:H^{1,p}(B_+)\to H^{1,p}(B_+)
\quad\hbox{and}\quad {\mathcal N}({\mathcal L}-{\mathcal L}_0)\phi_1:H^{1,2}(B_+)\to H^{1,2}(B_+)
\end{equation}
 have operator norms
 less than $1/2$. If we can show the right hand side of \eqref{I+S} is in $ H^{1,p}(B_+)$,
 then we can use a Neumann series to conclude  $\phi_0 u\in H^{1,p}(B_+)$. So we only need show
 \begin{equation}\label{2show}
 [{\mathcal L},\phi_0](\phi_1 u)\in H^{-1,p}(B_+).
 \end{equation}
 
 For $v\in H^{1,p'}(B_+)$, let us compute 
 \[
 \langle [{\mathcal L},\phi_0 ](\phi_1 u),v\rangle=
 \int_{B_+} a_{ij}(x)\left( \partial_j(\phi_1u)\partial_i(\phi_0 v)-\partial_j(\phi_0 u)\partial_iv\right)\,dx.
 \]
 Using $\phi_0\,\partial_ju\,\partial_iv=\phi_1\phi_0\,\partial_ju\,\partial_iv$, we find
 \[
 \left|  \langle [{\mathcal L},\phi_0 ](\phi_1 u),v\rangle \right|
 \leq C \int_{B_+} \left(|u|\,|\nabla v|+|\nabla u|\,|v|\right)\,dx.
 \]
 Let us first assume $n>2$. Then, by the H\"older and Sobolev inequalities,
 \[
  \int_{B_+} |u|\,|\nabla v|\,dx \leq \|u\|_{L^p(B_+)}\|v\|_{H^{1,p'}(B_+)} \leq C\|u\|_{H^{1,2}(B_+)}\|v\|_{H^{1,p'}(B_+)} 
 \]
 provided $p\leq 2n/(n-2)$. Similarly, we can use the H\"older and Sobolev inequalities to estimate
 \[
  \int_{B_+}  |\nabla u|\,|v|\,dx \leq \|u\|_{H^{1,2}(B_+)}\|v\|_{L^2(B_+)} \leq C \|u\|_{H^{1,2}(B_+)}\|v\|_{H^{1,p'}(B_+)}
 \]
 provided $2\leq np'/(n-p')$. But we can easily see that $p\leq 2n/(n-2)$ is equivalent to $2\leq np'/(n-p')$, so
 we have shown \eqref{2show} for $p=2(1+\alpha)$ where $\alpha = 2/(n-2)$. This is an improvement over $p=2$, and we can iterate it a finite number of times to conclude  \eqref{2show} for any $p>2$. If $n=2$, then the above argument works for any $2<p<\infty$.

\section{Derivation of the dynamical system}
\label{DerivationDynSys}

In this appendix we provide the details behind the derivation of the dynamical system \eqref{ourODEsystem:joint} for a given  solution
$u$ of the variational problem (\ref{eq:HalfspaceProblem1}).
Starting from \eqref{spectral-decomp}, we calculate
\[
\partial_j u = u_0'(r)\,\theta_j + (\widetilde{v}\,'(r)\cdot \widetilde{x})\,\theta_j+\widetilde{v}_j(r)+\partial_j w 
\quad\hbox{for}\ j=1,\dots,n-1
\]
and  \[
\partial_n u=u_0'(r)\,\theta_n + (\widetilde{v}\,'(r)\cdot \widetilde{x})\,\theta_n+\partial_n w.
\]
Now let us consider $\eta=\eta(r)$ in  (\ref{eq:HalfspaceProblem1}). Then $\partial_i\eta=\eta'(r)\,\theta_i$ for $i=1,\dots,n$, and plugging this and 
(\ref{spectral-decomp}) into (\ref{eq:HalfspaceProblem1}), we find
\begin{equation}\label{ode1}
\int_0^\infty \left[\left(\alpha\,u_0'+r\,\widetilde\beta\cdot\widetilde{v}\,'+\widetilde\gamma\cdot\widetilde v+p[\nabla w]-\overline{f}_1\right)\,\eta'+\overline{f_0}\,\eta\right]\,r^{n-1}\,dr=0,
\end{equation}
where
\[
\alpha(r)=\meanint_{S_+^{n-1}} \sum_{i,j=1}^n a_{ij}(r\theta)\theta_i\theta_j\,ds,
\quad 
\widetilde\beta_k(r)=\meanint_{S_+^{n-1}} \sum_{i,j=1}^n a_{ij}(r\theta)\theta_i\theta_j\theta_k\,ds \quad (k=1,\dots ,n-1),
\]
\[
\widetilde\gamma_j(r)=\meanint_{S_+^{n-1}}\sum_{i=1}^n a_{ij}(r\theta)\,\theta_i\, ds\quad (j=1,\dots ,n-1),
\quad
p[\nabla w](r)=\meanint_{S_+^{n-1}} \sum_{i,j=1}^n a_{ij}(r\theta)\,\partial_jw(r\theta)\,\theta_i\,ds,
\]
\begin{equation}\label{def:f1,f0}
\overline{f}_1(r)=\meanint_{S_+^{n-1}}\sum_{i=1}^n f_i(r\theta)\theta_i\,ds, \quad\hbox{and}\quad
\overline{f_0}(r)=\meanint_{S_+^{n-1}} f_0(r\theta)\,ds.
\end{equation}
Note that  $\alpha$, $p[\nabla w]$, $\overline{f}_1$, and $\overline{f_0}$ are scalar-valued while $\widetilde\beta$ and $\widetilde\gamma$ are $(n-1)$-vector-valued. For $0<r<1/4$ we have $\overline{f}_1(r)=0$, while
using (\ref{aij-deltaij}) and properties discussed in \cite{MM3}, we see that the others satisfy
\[
|\alpha(r)-1|, |\widetilde\beta(r)|, |\widetilde\gamma(r)| \leq\om(r) \quad \hbox{for}\ 0<r<1,
\]
\[
|p[\nabla w](r)|\leq \om(r) \meanint_{S_+^{n-1}} |\nabla w(r\theta)|\,ds  \quad \hbox{for}\ 0<r<1.
\]
Using (\ref{eq:aij_for_r>2}) and $u=0$ for $|x|>1$, we see that $\alpha(r)=1$ and $\widetilde\beta(r)=\widetilde\gamma(r)=p[\nabla w](r)=\overline{f}_1(r)=\overline{f_0}(r)=0$ for $r>1$.
Now if we integrate by parts in (\ref{ode1}) we obtain
\[
\int_0^\infty -[r^{n-1}(\alpha u_0'+r\widetilde\beta\cdot\widetilde v'+\widetilde\gamma\cdot\widetilde v+p[\nabla w]-\overline{f}_1)]'\eta
+r^{n-1}\,\overline{f_0}\,\eta \,dr=0,
\]
which means
\[
-[r^{n-1}(\alpha u_0'+r\widetilde\beta\cdot\widetilde v'+\widetilde\gamma\cdot\widetilde v+p[\nabla w]-\overline{f}_1)]'
+r^{n-1}\overline{f_0}=0.
\]
But we can integrate this to find
\begin{subequations}
\begin{equation}\label{eqnfor-u0'}
\alpha u_0'+r\widetilde\beta\cdot\widetilde v'+\widetilde\gamma\cdot\widetilde v+p[\nabla w]=\vartheta(r)
\end{equation}
where
\begin{equation}\label{def:vartheta}
\vartheta(r)=\overline{f}_1(r)+r^{1-n}\int_0^r\overline{f_0}(\rho)\rho^{n-1}\,d\rho.
\end{equation}
\end{subequations}
Since $\alpha(r)\geq \e>0$, \eqref{eqnfor-u0'} can be solved for $u_0'$ in terms of $\tilde v$ and $w$. 

Similarly, we can let $\eta=\eta(r)\,x_\ell$ in (\ref{eq:HalfspaceProblem1}) for $\ell=1,\dots,n-1$; this will give us a system of equations for the vector function $\widetilde v$. To begin with, we have $\partial_i\,\eta=r\eta'(r)\theta_i\theta_\ell+\eta(r)\delta_{i\ell}$. If we plug this and (\ref{spectral-decomp}) into (\ref{eq:HalfspaceProblem1}), we find
\[
\int_0^\infty \left[(u_0'\widetilde\beta+rA\widetilde v{\,}'+B\widetilde v+\widetilde\xi[\nabla w]-\widetilde f{\,}^\#)r\eta' +
(u_0'\widetilde\gamma+rB\widetilde v{\,}'+C\widetilde v+\widetilde\zeta[\nabla w]+\widetilde f{\,}^\flat)\eta\right] r^{n-1}dr=0,
\]
where $A$, $B$, and $C$ are $(n-1)\times(n-1)$ matrix-valued functions defined by
\begin{align}\notag
A_{\ell k}(r)&=\meanint_{S^{n-1}_+} a_{ij}(r\theta)\theta_i\theta_j\theta_\ell\theta_kds \quad (\ell,k=1,\dots ,n-1),\\ 
B_{\ell k}(r)&=\meanint_{S^{n-1}_+} a_{\ell j}(r\theta)\theta_j\theta_kds\quad (\ell,k=1,\dots ,n-1),\\ 
C_{\ell k}(r)&=\meanint_{S^{n-1}_+}  a_{\ell k}(r\theta)\,ds\quad (\ell,k=1,\dots ,n-1), \notag
\end{align}
and $\widetilde\xi[\nabla w]$, $\widetilde\zeta[\nabla w]$, $f^\#$, and $f^\flat$ are $(n-1)$-vector-valued functions defined by
\begin{align}\notag
\widetilde\xi_\ell[\nabla w](r)&=\meanint_{S^{n-1}_+} \sum_{i,j=1}^n a_{ij}\,\theta_i\,\theta_\ell\,\del_jw\,ds_\theta,
\quad
\widetilde\zeta_\ell[\nabla w](r)=\meanint_{S^{n-1}_+} \sum_{j=1}^n a_{\ell j}\,\del_jw \,ds,\notag \\
\widetilde f^\#_{\ell}(r)&=\meanint_{S^{n-1}_+} \sum_{i=1}^nf_i(r\theta)\theta_i\theta_\ell\,ds \quad\hbox{and}\quad
\widetilde f^\flat_{\ell}(r)=\meanint_{S^{n-1}_+}f_0(r\theta)\theta_\ell\,ds.
\end{align}
Using (\ref{aij-deltaij}) and properties discussed in \cite{MM3}, we see that these functions satisfy 
\begin{equation}
\begin{aligned}
\pmb{|}\,A-n^{-1}I_{n-1}\,\pmb{|},\ \pmb{|}\,B-n^{-1}I_{n-1}\,\pmb{|},\ \pmb{|}\,C-I_{n-1}\,\pmb{|} \leq \om(r)\quad \hbox{for}\ 0<r<1\\
|\widetilde\xi[\nabla w](r)|,\  |\widetilde\zeta[\nabla w](r)|
\leq \om(r)\,\meanint_{S_+^{n-1}} |\nabla w|\,ds\quad \hbox{for}\ 0<r<1.
\end{aligned}
\label{est:A-etc}
\end{equation}
(Here $I_{n-1}$ denotes the $(n-1)\times(n-1)$ identity matrix.)
For $r>1$ we use (\ref{eq:aij_for_r>2}) and $u=0$ to conclude
$A(r)=n^{-1}I=B(r)$, $C(r)=I$, and $\widetilde\xi[\nabla w](r)=0=\widetilde\zeta[\nabla w](r)=\widetilde f^\#(r)=\widetilde f^\flat(r)$. Now, using integration by parts, we obtain the 2nd-order system of ODEs
\begin{equation}
-\left[ r^n(u_0'\widetilde\beta+rA\widetilde v{\,}'+B\widetilde v+\widetilde\xi[\nabla w]-\widetilde f{\,}^\#)\right]'+
r^{n-1}(u_0'\widetilde\gamma+rB\widetilde v{\,}'+C\widetilde v+\widetilde\zeta[\nabla w]+\widetilde f{\,}^\flat)=0.
\label{eq:ODE1}
\end{equation} 

At this point we can use (\ref{eqnfor-u0'}) to eliminate $u_0'$ from (\ref{eq:ODE1}), and then use the change of variables $r=e^{-t}$.  After the change of variables we have
\[
\begin{aligned}
&\left[e^{-nt}\left(-A\widetilde v_t+B\widetilde v+\widetilde\xi[\nabla w]-\widetilde f{\,}^\#-\frac{1}{\alpha}(-\widetilde\beta\cdot\widetilde v_t+\widetilde\gamma\cdot\widetilde v+p[\nabla w]-\vartheta)\widetilde\beta\right)\right]_t \cr
+&\ e^{-nt}\left(-B\widetilde v_t+C\widetilde v+\widetilde\zeta[\nabla w]+\widetilde f{\,}^\flat-\frac{1}{\alpha}(-\widetilde\beta\cdot\widetilde v_t+\widetilde\gamma\cdot\widetilde v+p[\nabla w]-\vartheta)\widetilde\gamma\right)=0,
\end{aligned}
\]
which after some rearrangement can be written
\[
\begin{aligned}
\left[-A\widetilde v_t+B\widetilde v+\widetilde\xi[\nabla w]-\widetilde f{\,}^\#+\frac{\widetilde\beta\cdot\widetilde v_t-\widetilde\gamma\cdot\widetilde v-p[\nabla w]+\vartheta}{\alpha}\widetilde\beta\right]_t -(B-nA)\widetilde v_t+\frac{\widetilde\beta\cdot\widetilde v_t}{\alpha}(\widetilde\gamma-n\widetilde\beta)&
 \cr
+(C-nB)\widetilde v-\frac{\widetilde\gamma\cdot\widetilde v}{\alpha}(\widetilde\gamma-n\widetilde\beta)
=n\left[\widetilde\xi[\nabla w]-\frac{p[\nabla w]-\vartheta}{\alpha}\widetilde\beta\right]
+\frac{p[\nabla w]-\vartheta}{\alpha}\widetilde\gamma
-\widetilde\zeta[\nabla w]-\widetilde f{\,}^\flat&.
\end{aligned}
\]
To avoid differentiating the coefficient matrices, let us convert this to a first-order system for the $2(n-1)$-vector function $V=(V_1,V_2)$ where $V_1=\widetilde v$ and 
\[
V_2= -A\widetilde v_t+B\widetilde v+\widetilde\xi[\nabla w]-\widetilde f{\,}^\#+\frac{\widetilde\beta\cdot\widetilde v_t-\widetilde\gamma\cdot\widetilde v-p[\nabla w]+\vartheta}{\alpha}\widetilde\beta\,.
\]
Notice that the matrix $A$ is invertible near $x=0$ and for $|x|>1$, so we may assume that the variables were rescaled to make $A$ invertible for all $x$. Thus we may solve for the $t$-derivatives of $V_1$ and $V_2$ (which we now denote by the dot notation) to find:
\begin{equation}\label{eq:ODE2}
\begin{aligned}
\dot{V}_1 &-A^{-1}B V_1 +A^{-1}V_2-\frac{\widetilde\beta\cdot\dot{V}_1-\widetilde\gamma\cdot{V}_1}{\alpha}A^{-1}\widetilde\beta= A^{-1}\left[\widetilde\xi[\nabla w]-\widetilde f{\,}^\#-
\frac{p[\nabla w]-\vartheta}{\alpha}\widetilde\beta\right]\\
\dot{V}_2 &+(C-BA^{-1}B) V_1+(BA^{-1}-n) V_2+
\frac{\widetilde\beta\cdot\dot V_1}{\alpha}(\widetilde\gamma-(n+A^{-1})\vec\beta)
\\
&+\frac{\widetilde\gamma\cdot V_1}{\alpha}((n+A^{-1})\widetilde\beta-\widetilde\gamma)= n\left[\widetilde\xi[\nabla w]-\frac{p[\nabla w]-\vartheta}{\alpha}\widetilde\beta\right]
+\frac{p[\nabla w]-\vartheta}{\alpha}\widetilde\gamma
-\widetilde\zeta[\nabla w]-\widetilde f{\,}^\flat.
\end{aligned}
\end{equation}

Now \eqref{eq:ODE2} is still pretty complicated, but notice that the terms involving $\dot V_1$ and $\dot V_2$ in (\ref{eq:ODE2}) are of the form $(I+D(t))\dot V$ where $I$ is the $2(n-1)\times 2(n-1)$-identity matrix and the matrix $D(t)$ has matrix norm satisfying $\pmb{\bigr |}D(t)\pmb{\bigr |}\leq c\,\e^2(t)$. 
Consequently, we can multiply (\ref{eq:ODE2}) by $(I+D(t))^{-1}$ and,
after some calculations, see that $V$ satisfies a 1st-order system in the form
\begin{subequations}\label{ODE3:joint}
\begin{equation} \label{ODE3:a}
\frac{dV}{dt}+M(t)
V= F(t,\nabla w)+ F_0(t),
\end{equation}
where  $M(t)$ is a $2(n-1)\times 2(n-1)$ matrix of the form
\begin{equation}\label{ODE3:b}
\begin{aligned}
M(t)&=\, M_\infty +S_1(t)+S_2(t), \ \hbox{where}\\
M_\infty&=\begin{pmatrix}
 -I & n\,I \\ \frac{n-1}{n}\,I & (1-n)\,I
 \end{pmatrix} \ \hbox{and}\\
  S_1(t)&=
\begin{pmatrix}
I-A^{-1}B & A^{-1}-nI \\ C-BA^{-1}B+\frac{1-n}{n}I & BA^{-1}-I 
\end{pmatrix}.
\end{aligned}
 \end{equation}
The $S_i$ satisfy
 \begin{equation}\label{ODE3:c}
 \begin{aligned}
\pmb{\bigr |}S_1(t)\pmb{\bigr |}\leq\e(t) \  \hbox{and}\ \pmb{\bigr |}S_2(t)\pmb{\bigr |}\leq c\,\e^2(t) \ \hbox{as}\ t\to\infty,\\
 S_1(t)=0=S_2(t)  \ \hbox{for}\ t<0,
 \end{aligned}
 \end{equation}
while the vector  $F(t,\nabla w)$ satisfies
 \begin{equation}\label{ODE3:d}
 | F(t,\nabla w)|\leq \ c\, \e(t)\,\meanint_{S_+^{n-1}} |\nabla w|\,ds\ \hbox{as}\ t\to\infty
\  \ \hbox{and}\ \ 
F(t,\nabla w)\equiv 0\  \hbox{for}\ t<0,
 \end{equation}
 and the vector $F_0(t)$ has support in $t>0$ with $L^1$-norm satisfying
  \begin{equation}\label{ODE3:e}
 \| F_0\|_{L^1(\RR)}\leq \ c\, (\|\vec f\,\|_p+\|f_0\|_p).
 \end{equation}
 \end{subequations}
 Note that $M(t)$ and $F_0(t)$ depend on $a_{ij}$, $\vec f$, and $f_0$, but not on $w$.
 
 We can further simplify our dynamical systems by another change of dependent variables. We can calculate the 
eigenvalues of $M_\infty$ to be $\lambda=0$ and $\lambda=-n$ (each occurring $n-1$ times). The matrix
\[
J=\begin{pmatrix} nI & nI \\ I  & (1-n)I \end{pmatrix}
\]
diagonalizes $M_\infty$, i.e.\ $J^{-1}M_\infty J=\hbox{diag}(0,\dots,0,-n,\dots,n)$, so let us introduce new dependent variables $V\to(\phi,\psi)$ by
\begin{equation}
V=J\begin{pmatrix} \phi \\ \psi \end{pmatrix}.
\end{equation}
We find that the dynamical system (\ref{ODE3:a}) now takes the form
  (\ref{ourODEsystem:a}), where the conditions (\ref{ourODEsystem:c}) and  (\ref{ourODEsystem:d}) follow from 
  (\ref{ODE3:d}) and (\ref{ODE3:e}) respectively, 
  and ${\mathcal R}$ is of the form (\ref{ourODEsystem:b}) 
  with
 \begin{equation}
R_1=\frac{n-1}{n^2}A^{-1}-\frac{n-1}{n}A^{-1}B+C-BA^{-1}B+\frac{1}{n}BA^{-1}-I.
\end{equation}
To simplify this expression for $R_1$, let us write
\[
A=n^{-1}(1+\widetilde A), 
\qquad
B=n^{-1}(1+\widetilde B), 
\qquad\hbox{and}\quad
C=n^{-1}(1+\widetilde C),
\]
where
$\pmb{|}\widetilde A\pmb{|}, \  \pmb{|}\widetilde B\pmb{|}, \  \pmb{|}\widetilde C\pmb{|} \leq c\,\e(t)$ as $t\to\infty$.
Then
$
A^{-1}\approx n(I-\widetilde A),
$
 and a calculation shows
\[
R_1\approx \widetilde C-\widetilde B=C-nB\quad\hbox{as}\ t\to\infty,
\]
which gives the formula \eqref{ourODEsystem:c}. Finally, if we follow our changes of dependent variables from $(\widetilde v,\widetilde v_r)$ to $(\varphi,\psi)$, we easily see that \eqref{ourODEsystem:f} holds.

Now our original assumption that $u\in H^{1,2}_\loc(\RR^n_+)$ has implications for $V(t)$ as $t\to\infty$. In fact, using orthogonality properties in the decomposition \eqref{spectral-decomp}, we find that $\nabla u\in L^2(B_+^n)$, where $B_+^n=B_1(0)\cap\RR^n_+$, implies
\[
\int_0^1 \left( (u_0')^2+|\widetilde v|^2+r^2|\widetilde v\,'|^2\right)r^{n-1}\,dr<\infty
\quad\hbox{and}\quad
\nabla w\in L^2(B_+^n).
\]
In particular, using $V_1=\widetilde v$ and the second equation in \eqref{eq:ODE2} for $\dot{V_2}$, this implies
\begin{equation}\label{eq:L2-boundonV}
\int_0^\infty \left( |V|^2+|\dot{V}|^2+|\overline{\nabla w}|^2\right)e^{-nt}\,dt <\infty,
\end{equation}
where 
\[
\overline{\nabla w}=\meanint_{S^{n-1}_+}\nabla w\,ds.
\]
Thus $V(t)$ and its first-order derivative cannot grow too rapidly as $t\to\infty$.

\section{Stability properties of dynamical systems}\label{StabilityPropsDynSys}

Here we recall a result on stability properties of dynamical systems that was obtained in \cite{MM3}. 
Let $\e(t)$ be a positive, nonincreasing continuous function satisfying
 \[
 \int_0^\infty \e^2(\tau)\,d\tau<\infty.
 \]
Consider a $2k\times 2k$-dimensional dynamical system in the form
 \begin{subequations} \label{DynSys}
 \begin{equation}  \label{DynSys:a}
\frac{d}{dt}
\begin{pmatrix}
\varphi \\ \psi
\end{pmatrix}
+
\begin{pmatrix}
0 & 0 \\ 0 & -nI
\end{pmatrix}
\begin{pmatrix}
\varphi \\ \psi
\end{pmatrix}
+
{\mathcal R}(t)
\begin{pmatrix}
\varphi \\ \psi
\end{pmatrix}
=  
g(t) \quad\hbox{for}\ t>0,
\end{equation}
where $n>0$ and ${\mathcal R}$ can be written as a matrix of $k\times k$-blocks 
\begin{equation}\label{DynSys:b}
{\mathcal R}(t)=
\begin{pmatrix}
R_{1}(t) & R_{2}(t) \\ R_{3}(t) & R_{4}(t)
\end{pmatrix}
\quad
\hbox{with}\ \pmb{\bigr |}R_j(t)\pmb{\bigr |}\leq \e(t)\ \hbox{on}\ 0<t<\infty,
\end{equation}
with the block $R_1$ satisfying 
\begin{equation}\label{DynSys:c}
|R_1(t)-R(t)|\leq c\,\e^2(t) \ \hbox{as}\ t\to\infty, 
\end{equation}
for a certain $k\times k$ matrix $R(t)$. We also assume that the vector function $g(t)=(g_1(t),g_2(t))$ satisfies the following conditions:
\begin{equation} \label{g1-condition}
g_1\in L^1(0,\infty)
\end{equation}
and there exists $\de>0$ such that for any choice of $\alpha\in [n-\de,n)$ there is a constant $c_\alpha$ so that
\begin{equation}\label{g2-condition}
e^{\alpha t}\int_t^\infty |g_2(s)|\,e^{-\alpha s}\,ds \leq c_\alpha \e(t) \quad\hbox{for}\ 0<t<\infty.
\end{equation}
\end{subequations}
We want to relate the stability for \eqref{DynSys} to that for 
\begin{subequations}\label{DynSys2}
\begin{equation}\label{DynSys2:a}
\frac{d\varphi}{dt}+R\varphi=0 \quad\hbox{for}\ t>0,
\end{equation}
and the ``finite-energy'' condition on $\psi$
\begin{equation}\label{DynSys2:b}
\int_0^\infty \left(|\psi|^2+|\psi_t|^2\right)\,e^{-nt}\,dt <\infty.
\end{equation}
\end{subequations}

\begin{proposition}\label{pr:stabilty}  
Suppose that ${\mathcal R}$ and $g=(g_1,g_2)$ satisfy \eqref{g1-condition} and \eqref{g2-condition}. Assume also that \eqref{DynSys2:a} is uniformly stable. Then all solutions $(\phi,\psi)$ of \eqref{DynSys:a} that satisfy \eqref{DynSys2:b} will remain bounded as $t\to\infty$, and $\psi(t)\to 0$. In fact, for $\alpha=n-\delta$ with $\delta>0$ sufficiently small, we will have the estimates
\begin{subequations}
\begin{equation}\label{est:phi}
\sup_{0<t<\infty} |\varphi(t)|\leq c\,(c_\alpha+|\varphi(0)|+\|g_1\|_1),
\end{equation}
\begin{equation}\label{est:psi}
 |\psi(t)|\leq c\,\e(t)(c_\alpha+\sup_{t<\tau<\infty}|\varphi(\tau)|).
\end{equation}
\end{subequations}
In addition, if all solutions of \eqref{DynSys2:a} are asymptotically constant as $t\to\infty$, then the solution $(\varphi,\psi)$ of 
\eqref{DynSys} also has a limit:
\begin{equation}
(\varphi(t),\psi(t))\to (\varphi_\infty,0) \quad\hbox{as}\ t\to\infty.
\end{equation}
\end{proposition}

\begin{remark} \label{rk:1} 
In \cite{MM3}, Proposition \ref{pr:stabilty} was stated and proved for the special case $k=n$. However, the proof in \cite{MM3} does not use $k=n$, so so it proves the above Proposition. This is important for the application in this paper since we need to take $k=n-1$. 
\end{remark}

\begin{remark} \label{rk:2} 
Since our dynamical system \eqref{DynSys2:a} is linear, the condition that it be uniformly stable is equivalent to the condition $|\Phi(t)\Phi^{-1}(s)|\leq K$ for $t>s>0$, where $\Phi$ denotes the fundamental matrix for \eqref{DynSys2:a}. {\rm (Cf.\ \cite{C}.)}
\end{remark}

\section{Uniqueness of solutions}\label{Uniqueness}

In this appendix we discuss uniqueness for solutions of our variational equation. 
Suppose $u\in H^{1,p}_\loc(\overline{\RR^n_+})$ for $p\geq 2$ satisfies
\begin{equation}\label{weaksolution-homog}
\int_{\RR^n_+} a_{ij}\partial_j u\, \partial_i\eta \, dx=0 \quad\hbox{for}\ \eta\in C^1_{comp}(\overline{\RR^n_+}),
\end{equation}
and 
\begin{equation}\label{M1p-est}
M_{1,p}(u,r)\leq C\,r^{-\alpha} \quad\hbox{for}\ r>1.
\end{equation}
The following proposition describes that values of $\alpha>0$ for which we can conclude that $u\equiv 0$.
\begin{proposition} \label{pr:uniqueness} 
Suppose $u\in H^{1,p}_\loc(\overline{\RR^n_+})$ for $p\geq 2$ satisfies
\eqref{weaksolution-homog} and \eqref{M1p-est} where $\alpha$ satisfies
\begin{equation}\label{alpha>}
\alpha > \frac{n(p-2)}{2p}.
\end{equation}
Then $u\equiv 0$.
\end{proposition}
\noindent
As a special case we obtain the uniqueness result that is useful in our proof of Theorem \ref{th:1}.
\begin{corollary}\label{co:uniqueness} 
If $u\in H^{1,p}_{\loc}(\RR^n_+)$  for $p\geq 2$ is a solution of \eqref{eq:HalfspaceProblem1} that satisfies
\[
M_{1,p}(u,r)\leq C\,r^{-n} \quad\hbox{for}\ r>1,
\]
then $u$ is unique.
\end{corollary}

\noindent
{\bf Proof of Proposition \ref{pr:uniqueness}.} The strategy is to show that \eqref{weaksolution-homog} holds with $\eta=u$, i.e.\
\begin{equation}\label{u=eta}
\int_{\RR^n_+} a_{ij}\partial_j u\, \partial_i u \,dx=0.
\end{equation}
The ellipticity of $a_{ij}$ then implies $\nabla u\equiv 0$, i.e.\ that $u$ is constant. Finally, $\alpha>0$ in \eqref{M1p-est} implies that $u\equiv 0$.

First let us determine the values of $\alpha>0$ that imply that $\nabla u\in L^2(\RR^n_+)$. Since $p\geq 2$, we know that $\nabla u\in L^2_{\loc}(\overline{\RR^n_+})$, so the question is whether
\begin{equation}\label{u-finite_energy}
\int_{x\in\RR^n_+,\, |x|>1} |\nabla u|^2\,dx <\infty.
\end{equation}
But 
\[
\begin{aligned}
\int_{x\in\RR^n_+,\, |x|>1} |\nabla u|^2\,dx = \int_1^\infty \int_{S^{n-1}_+} |\nabla u(r\theta)|^2 \,r^{n-1}\,ds\,dr \\
= \sum_{j=0}^\infty \int_{2^j}^{2^{j+1}} \int_{S^{n-1}_+} |\nabla u|^2 \,r^{n-1}\,ds\,dr,
\end{aligned}
\]
and by H\"older's inequality
\[
\begin{aligned}
\int_{2^j}^{2^{j+1}} & \int_{S^{n-1}_+} |\nabla u|^2 \,r^{n-1}\,ds\,dr=\int_{2^j<|x|<2^{j+1}}|\nabla u|^2\,dx \\
& \leq \left( \int_{2^j<|x|<2^{j+1}} |\nabla u|^p\,dx\right)^{2/p} \left( \int_{2^j<|x|<2^{j+1}} dx\right)^{(p-2)/p} \\
& \leq |S^{n-1}_+| \,M_{1,p}(u,2^j)^2 \, 2^{nj(p-2)/p} \ \leq \ C\,2^{j(-2\alpha+\frac{n(p-2)}{p})}.
\end{aligned}
\]
Thus, we see that \eqref{alpha>} implies \eqref{u-finite_energy}.

Now let us verify that \eqref{u-finite_energy} is sufficient to enable us to take $u=\eta$ in \eqref{weaksolution-homog}, i.e.\ that \eqref{u=eta} holds. It suffices to show that there exist $u_m\in C^1_{comp}(\overline{\RR^n_+})$ with $\nabla u_m\to\nabla u$ in $L^2(\RR^n_+)$ as $m\to \infty$. But, using mollifiers, it suffices to show that this can be achieved with $u_m\in H^{1,2}_{comp}(\overline{\RR^n_+})$. So let $\chi(t)$ be a smooth function for $t>0$ with $|\chi'(t)|\leq 2$ and
\[
\chi(t)=\begin{cases} 0 & \hbox{if}\  t>2 \\ 1 & \hbox{if}\  0<t<1.
\end{cases}
\]
Then, for $m=1,2,\dots$, define $u_m\in H^{1,2}_{comp}(\overline{\RR^n_+})$ by
\[
u_m(x)=u(x)\cdot \chi_m(|x|) \quad\hbox{where}\ \chi_m(t)=\chi(t/m).
\]
For $i=1,\dots,n$ we compute
\[
\nabla_i u(x)-\nabla_i u_m(x) = (1-\chi_m(|x|))\nabla_i u(x) + u(x)\cdot\chi_m'(|x|)\cdot\frac{x_i}{|x|}.
\]
We want to show both terms on the right tend to zero in $L^2(\RR^n_+)$ as $m\to\infty$.
If we assume \eqref{alpha>}, then we know $\nabla u\in L^2(\RR^n_+)$, and hence
\[
\int_{\RR^n_+} (1-\chi_m)^2|\nabla u|^2\,dx\leq \int_{\RR^n_+,|x|>m}|\nabla u|^2\,dx \to 0 \quad\hbox{as}\ m\to\infty.
\]
To estimate the second term we use
\[
\begin{aligned}
\int_{\RR^n_+} (\chi'_m)^2|u|^2\,dx &\leq \frac{4}{m^2}\int_{m<|x|<2m}|u(x)|^2\,dx \\
&\leq \frac{C}{m^2}\left(\int_{m<|x|<2m}|u|^p\,dx\right)^{2/p}\left(m^n\right)^{(p-2)/p} \\
&\leq C\, m^{n-2-\frac{2n}{p}} M_{1,p}(u,m)^2,
\end{aligned}
\]
which tends to zero as $m\to\infty$ provided $n-2-(2n/p)-2\alpha<0$, i.e.\
\[
\alpha> \frac{n(p-2)}{2p} -1.
\]
But this condition on $\alpha$ is certainly implied by \eqref{alpha>}, so we are done. $\Box$


\end{document}